\documentclass[11pt]{article}
\usepackage{amsmath,amssymb,latexsym,graphics,epsfig,psfrag,float}
\begin{document}

\begin{center}
{\Large{\bf Augmented Sparse Principal Component Analysis for High Dimensional
Data}}\footnote{This manuscript was written in 2007, and a version dated
December 11, 2007 has been available on the first author's website at
\texttt{http://anson.ucdavis.edu/$\sim$debashis/techrep/augented-spca.pdf}. But
it is posted to arXiv now in its 2007 form as that version has since been cited
in work by many authors. Revisions incorporating later work will be posted
separately.}

\vskip.1in Debashis Paul$^\dag$ and Iain M. Johnstone$^\ddag$

\vskip.1in \dag \textit{University of California, Davis}

\ddag \textit{Stanford University}

\end{center}

\begin{abstract}
We study the problem of estimating the leading eigenvectors of a
high-dimensional population covariance matrix based on independent Gaussian
observations. We establish lower bounds on the rates of convergence of the
estimators of the leading eigenvectors under $l^q$-sparsity constraints when an
$l^2$ loss function is used. We also propose an estimator of the leading
eigenvectors based on a coordinate selection scheme combined with PCA and show
that the proposed estimator achieves the optimal rate of convergence under a
sparsity regime. Moreover, we establish that under certain scenarios, the usual
PCA achieves the minimax convergence rate.
\end{abstract}

\section{Introduction}

Principal components analysis (PCA) has been a widely used technique
in reducing dimensionality of multivariate data. A traditional
setting where PCA is applicable is when one has repeated
observations from a multivariate population that can be described
reasonably well by its first two moments. When the dimension of
sample observations, is fixed, distributional properties of the
eigenvalues and eigenvectors of the sample covariance have been
dealt with at length by various authors. Anderson (1963), Muirhead
(1982) and Tyler (1983) are among standard references. Much of the
``large sample'' study of the eigen-structure of the sample
covariance matrix is based on the fact that, sample covariance
approximates population covariance matrix well when sample size is
large. However, due to advances in data acquisition technologies,
statistical problems, where the dimensionality of individuals are of
nearly the same order of magnitude as (or even bigger than) the
sample size, are increasingly common. The following is a
representative list of areas and articles where PCA has been in use.
In all these cases $N$ denotes the dimension of an observation and
$n$ denotes the sample size.
\begin{itemize}
\item \textit{Image recognition :} The face recognition problem is to
  identify a face from a collection of faces. Here each observation is
  a digitized image of the face of a person. So typically, with
  $128\times 128$ pixel grids, one has to deal with a situation
  where $N \approx 1.6 \times 10^{6}$. Whereas, a standard image
  database, e.g. that of students of Brown University
  Wickerhauser (1994), may contain only a few hundred pictures.
\item \textit{Shape analysis :} Stegmann and Gomez (2002),
  Cootes, Edwards and Taylor (2001) outline a
  class of methods for analyzing the shape of an object based on
  repeated measurements that involves annotating the objects for
  landmarks. These landmarks act as features of the objects,
  and hence, can be thought of as the dimension of the observations.
  For a specific example relating to motion of hand Stegmann and Gomez (2002),
  the number of landmarks
  is 56 and sample size is 40.
\item \textit{Chemometrics :} In many chemometric studies, sometimes
  the data consists of several thousands of spectra measured at
  several hundred wavelength positions, e.g. data collected
  for calibration of
  spectrometers. Vogt, Dable, Cramer and Booksh (2004) give an overview of some of
  these applications.
\item \textit{Econometrics :} Large factor analysis models are often
  used in econometric studies, e.g. in dealing with hundreds of stock
  prices as a multivariate time series. Markowitz's theory of optimal
  portfolios ask this question. \textit{Given a set of financial assets
  characterized by their average return and risk, what is the optimal
  weight of each asset, such that the overall portfolio provides the
  best return?} Laloux, Cizeau, Bouchaud and Potters (2000) discuss
  several applications.
  Bai (2003) considers some inferential aspects.
\item \textit{Climate studies :} Measurements on atmospheric
  indicators, like ozone concentration etc. are taken at a number
  of monitoring stations over a number of time points. In this
  literature, principal
  components are commonly referred to as ``empirical orthogonal
  functions''. Preisendorfer (1988) gives a detailed
  treatment. EOFs are also used for model diagnostics and data summary
  Cassou, Deser, Terraty, Hurrell and Dr\'{e}villon (2004).
\item \textit{Communication theory :} Tulino and Verdu (2004)
  give an extensive treatment to the connection between random
  matrix theory and vector channels used in wireless communications.
\item \textit{Functional data analysis :} Since observations are
  curves, which are typically measured at a large number of points, the
  data is high dimensional.  Buja, Hastie and Tibshirani (1995) give an example of
  speech dataset consisting of 162 observations - each one is a
  periodogram of a ``phoneme'' spoken by a
  person. Ramsay and Silverman (2002) discuss other applications.
\item \textit{Microarray analysis :} Gene microarrays present data in
  the form expression profiles of several thousand genes for each
  subject under study. Bair, Hastie, Paul and Tibshirani (2006) analyze an example
  involving the study of survival times of 240 ($= n$) patients
  with diffuse large B-cell lymphoma, with gene expression measurements
  for 7389 ($= N$) genes.
\end{itemize}
Of late, researchers in various fields have been using different
versions of non-identity covariance matrices of growing dimension.
Among these, a particularly interesting model assumes that,
\begin{itemize}
\item[(*)] the eigenvalues of the population covariance matrix $\Sigma$
are (in descending order)
$$
\ell_1, \ldots, \ell_M, \sigma^2, \ldots,\sigma^2,
$$
where $\ell_M > \sigma^2 > 0$.
\end{itemize}
This has been deemed the ``spiked population model'' by Johnstone
(2001). It has also been observed that for certain types of data,
e.g. in speech recognition Buja, Hastie and Tibshirani (1995),
wireless communication Telatar (1999), statistical learning (Hoyle
and Rattray (2003, 2004)), a few of the sample eigenvalues have
limiting behavior that is different from the behavior when the
covariance is the identity. This paper deals with the issue of
estimating the eigenvectors of $\Sigma$, when it has the structure
described by (*), and the dimension $N$ grows to infinity together
with sample size $n$.

In many practical problems, at least the leading eigenvectors are
thought to represent some underlying phenomena. This has been one of
the reasons for their popularity in analysis of what can be
characterized as functional data. For example, Zhao, Marron and
Wells (2004) consider the ``yeast cell cycle'' data of Spellman
\textit{et al.} (1998), and argue that the first two components
obtained by a functional PCA of the data represent systematic
structure. In climate studies, empirical orthogonal functions are
often used for identifying patterns in the data, as well as for data
summary. See for example Corti, Molteni and Palmer (1999). In many
of these instances there is some idea about the structure of the
eigenvectors of the covariance matrix, such as to the extent they
are smooth, or oscillatory. At the same time, these data are often
corrupted with a substantial amount of noise, which can lead to very
noisy estimates of the eigen-elements. There is also a growing
literature on functional response models in which the regressors are
random functions and the responses are either vectors or functions
(Chiou, M\"{u}ller and Wang (2004), Hall and Horowitz (2004),
Cardot, Ferraty and Sarda (2003)). Quite often a functional
principal component regression is used to solve these problems.
Thus, there are both practical and scientific interests in devising
methods for estimating the eigenvectors and eigenvalues that can
take advantage of the information about the structure of the
population eigenvectors. At the same time, there is also a need to
address this estimation problem from a broader statistical
perspective.

In multivariate analysis, there is a huge body of work on estimation
of population covariance, and in particular on developing optimal
strategies for estimation from a decision theoretic point of view.
Dey and Srinivasan (1985), Efron and Morris (1976), Haff (1980), Loh
(1988) are some of the standard references in this field. However, a
decision theoretic treatment of functional data analysis is still
somewhat limited in its breadth. Hall and Horowitz (2004) and Tony
Cai and Hall (2005) derive optimal rates of convergence of
estimators of the regression function and fitted response in
functional linear model context. Cardot (2000) gave upper bounds on
the rate of convergence of a spline-based estimator of eigenvectors
under some smoothness assumptions. Kneip (1994) also derived similar
results in a slightly different context.

In this paper, the aim is to address the problem of estimating
eigenvectors from a minimax risk analysis viewpoint.  Henceforth,
the observations will be assumed to have a Gaussian distribution.
This assumption, though somewhat idealized, helps in bringing out
some essential features of the estimation problem. Since algebraic
manipulation of spectral elements of a matrix is rather difficult,
it is not easy to make any precise finite sample statement about the
risk properties of estimators. Therefore the analysis is mostly
asymptotic in nature, even though efforts have been made to make the
approximations to risk etc. as explicit as possible. The asymptotic
regime considered here assumes a triangular array structure in which
$N$, the dimensionality of individual observations, tends to
$\infty$ with sample size $n$. This framework is partly motivated by
similar analytical approaches to the problem of estimation of mean
function in nonparametric regression context. In particular, a
squared error type loss is proposed, and some $l^q$-type sparsity
constraint is imposed on the parameters, which in our case are
individual eigenvectors. Relevance of this sort of constraints in
the  context of functional data analysis is discussed in Section
\ref{sec:highdpca-sparsity}. The main results of this chapter are
the following. Theorem 1
describes risk behavior of sample eigenvectors as estimators of
their population counterparts. Theorem 2
gives a lower bound on the minimax risk.  An estimation scheme,
named \textit{Augmented Sparse Principal Component Analysis (ASPCA)}
is proposed and is shown to have the optimal rate of convergence
over a class of $l^q$ norm-constrained parameter spaces under
suitable regularity conditions. Throughout it is assumed that the
leading eigenvalues of the population covariance matrix are
distinct, so the eigenvectors are identifiable. A more general
framework, which looks at estimating the eigen-subspaces and allows
for eigenvalues with arbitrary multiplicity, is beyond the scope of
this paper.

\section{Model}\label{sec:model}

Suppose that, $\{X_i : i=1,\ldots,n\}_{n \geq 1}$ is a triangular
array, where the $N \times 1$ vectors $X_i := X_{i}^n, i=1,\ldots,n$
are i.i.d. on a common probability space for each $n$. The dimension
$N$ is assumed to be a function of $n$ and increases without bound
as $n \to \infty$. The observation vectors are assumed to be i.i.d.
as $N(\xi,\Sigma)$, where $\xi$ is the mean vector; and $\Sigma$ is
the covariance matrix. The assumption on $\Sigma$ is that, it is a
finite rank perturbation of (a multiple of) the identity. In other
words,
\begin{equation}\label{eq:Sigma_basic}
\Sigma = \sum_{\nu=1}^M \lambda_\nu \theta_\nu \theta_\nu^T +
\sigma^2 I,
\end{equation}
where $\lambda_1 > \lambda_2 > \ldots > \lambda_M > 0$, and the
vectors $\theta_1,\ldots,\theta_M$ are orthonormal. Notice that
strict inequality in the order relationship among the
$\lambda_\nu$'s implies that the $\theta_\nu$ are identifiable up to
a sign convention. Notice that with this identifiability condition,
$\theta_\nu$ is the eigenvector corresponding to the $\nu$-th
largest eigenvalue, namely, $\lambda_\nu + \sigma^2$, of $\Sigma$.
The term ``finite rank'' means that, $M$ will remain fixed for all
the asymptotic analysis that follows. This analysis involves letting
both $n$ and $N$ increase to infinity simultaneously.  Therefore,
$\Sigma$, the $\lambda_\nu$'s and the $\theta_\nu$'s should be
thought of as being dependent on $N$.

The observations can be equivalently described in terms of the
\textit{factor analysis model} :
\begin{equation}\label{eq:basic}
X_{ik} = \xi + \sum_{k=1}^M \sqrt{\lambda_\nu} v_{\nu i} \theta_{\nu
k} + \sigma Z_{ik}, \quad i=1,\ldots,n, \quad k=1,\ldots,N.
\end{equation}
Here, for each $n$, $v_{\nu i}$, $Z_{ik}$ are all independently and
identically distributed as $N(0,1)$. $M \geq 1$ is assumed fixed.

Since the eigenvectors of $\Sigma$ are invariant to a scale change
in the original observations, for simplifying notation, it is
assumed that $\sigma = 1$. Notice that this also means that,
$\lambda_1,\ldots,\lambda_M$ appearing in the results relating to
the rates of convergence of various estimators of $\theta_\nu$
should be changed to $\lambda_1/\sigma, \ldots, \lambda_M/\sigma$
when (\ref{eq:Sigma_basic}) holds with an arbitrary $\sigma > 0$.

Another simplifying assumption is that, $\xi = 0$. This is because,
the main focus of the current exposition is on estimating the
eigen-structure of $\Sigma$, and the unnormalized  sample covariance
matrix
$$
\sum_{i=1}^n (X_i - \overline{X})(X_i - \overline{X})^T,
$$
where $\overline{X}$ is the sample mean, has the same distribution
as that of the matrix
$$
\sum_{i=1}^{n-1} Y_iY_i^T,
$$
where $Y_i$ are i.i.d. $N(0,\Sigma)$. This means that, for
estimation purposes, if the attention is restricted to the sample
covariance matrix, then from an asymptotic analysis point of view,
it is enough to assume $\xi = 0$, and to define the sample
covariance matrix as $\mathbf{S} =
\frac{1}{n}\mathbf{X}\mathbf{X}^T$, where $\mathbf{X} =
[X_1:\ldots:X_n]$.

The following condition, or \textit{Basic Assumption} will be used
frequently, and will be referred to as {\bf BA}.
\begin{itemize}
\item[{\bf BA}]
(\ref{eq:basic}) and (\ref{eq:Sigma_basic}) hold, with $\xi =0$ and
$\sigma = 1$; $N = N(n) \to \infty$ as $n \to \infty$; $\lambda_1 >
\ldots > \lambda_M > 0$.
\end{itemize}
For the estimation problem, it may be assumed that, as $n, N \to
\infty$, $\theta_\nu := \theta_\nu^n \to \overline\theta_\nu$ in
$l^2(\mathbb{R})$, though it is not strictly necessary. But this
assumption is appropriate if the observation vectors are the vectors
of first $N$ coefficients of some noisy function in $L^2(D)$ (where
$D$ is an interval in $\mathbb{R}$), when represented in a suitable
orthogonal basis for the $L^2(D)$ space. See Section
\ref{subsec:highdpca-fda_model} for more details. In such cases one
can talk about estimating the eigenfunctions of the underlying
covariance operator, and the term consistency has its usual
interpretation. However, even if $\theta_\nu^n$ does not converge in
$l^2$, one can still use the term ``consistency'' of an estimator
$\widehat \theta_\nu^n$ to mean that $L(\widehat \theta_\nu^n,
\theta_\nu^n) \to 0$ in probability as $n \to \infty$, where $L$ is
an appropriate loss function.

\subsection{Squared error type loss}\label{subsec:highdpca-lossfn}

The goal is, given data $X_1,X_2,\ldots,X_n$,  to estimate
$\theta_\nu$, for $\nu=1,\ldots,M$. To assess the performance of any
such estimator, a minimax risk analysis approach is proposed. The
first task is to specify a loss function for this estimation
problem. Observe that since the model is invariant under separate
changes of sign of the $\theta_\nu$, it is necessary to specify a
loss function that is also invariant under a sign change. We specify
the following loss function :
\begin{equation}\label{eq:loss_fn}
L(\mathbf{a},\mathbf{b})  = L([\mathbf{a}],[\mathbf{b}]) := 2(1 -
|\langle \mathbf{a},\mathbf{b}\rangle|) = \parallel \mathbf{a} -
sign(\langle \mathbf{a}, \mathbf{b} \rangle) \mathbf{b} \parallel^2,
\end{equation}
where $\mathbf{a}$ and $\mathbf{b}$ are $N \times 1$ vectors with
$l^2$ norm 1; and $[\mathbf{a}]$ denotes the equivalence class of
$\mathbf{a}$ under sign change. Note that,
$L(\mathbf{a},\mathbf{b})$ can also be written as $\min\{\parallel
\mathbf{a} - \mathbf{b}
\parallel^2, \parallel \mathbf{a}
  + \mathbf{b} \parallel^2\}$. There is another useful relationship with a
different loss function, denoted by $L_s(\mathbf{a},\mathbf{b}) :=
\sin^2 \angle (\mathbf{a},\mathbf{b})$, for any two  $N \times 1$
unit vectors $\mathbf{a}$ and $\mathbf{b}$. $\sin \angle
(\cdot,\cdot)$ is a metric on the space $\mathbb{S}^{N-1}$, i.e. the
unit sphere in $\mathbb{R}^N$. Also, $L_s(\mathbf{a},\mathbf{b}) =
\sin^2\angle (\mathbf{a},\mathbf{b}) = 1 - |\langle \mathbf{a},
\mathbf{b} \rangle|^2 =
L(\mathbf{a},\mathbf{b})(2-L(\mathbf{a},\mathbf{b}))$. Hence, if
$L(\mathbf{a},\mathbf{b}) \approx 0$, then these two quantities have
approximately the same value. This implies that, the asymptotic risk
bounds derived in terms of the loss function $L$ remain valid, up to
a constant factor, for the loss function $L_s$ as well.

\subsection{Rate of convergence for ordinary PCA}

It is assumed that either $\lambda_1$ is fixed, or that it varies
with $n$ and $N$ so that,
\begin{itemize}
\item[{\bf L1}]
as $n, N \to \infty$, $\frac{\lambda_\nu}{\lambda_1} \to \rho_\nu$
for $\nu=1,\ldots,M$, where $1 = \rho_1 > \rho_2 > \ldots > \rho_M$;
\item[{\bf L2}] as $n,N \to \infty$, $\frac{N}{nh(\lambda_1)} \to 0$,
  where
\begin{equation}\label{eq:h_lambda}
h(\lambda) = \frac{\lambda^2}{1+\lambda}~.
\end{equation}
\end{itemize}
Notice that, all four conditions (i)-(iv) below imply that
$\frac{N}{nh(\lambda_1)} \to 0$ as $n \to \infty$.
\begin{itemize}
\item[(i)]{$\frac{N}{n} \to \gamma \in (0,\infty)$ and
$\frac{N}{n\lambda_1} \to 0$}
\item[(ii)]{$\lambda_1 \to 0$, $\frac{N}{n} \to 0$ and
$\frac{N}{n \lambda_1^2} \to 0$}
\item[(iii)]{$0 < \lim\inf_{n \to \infty} \lambda_1 \leq
    \lim\sup_{n \to \infty} \lambda_1 < \infty$ and $\frac{N}{n} \to 0$}
\item[(iv)]{$\frac{N}{n} \to \infty$, and
$\frac{N}{n\lambda_1} \to 0$.}
\end{itemize}
\vskip.1in\noindent{\bf Remark :}\label{rem:cond_on_lambda}
Condition {\bf L1} is really an asymptotic identifiability condition
which guarantees that at the scale of the largest ``signal''
eigenvalue, bigger eigenvalues are well-separated.

\vskip,1in\noindent{\bf Theorem 1:} \label{thm:OPCA_risk_bound}
\textit{Suppose that the eigenvalues $\lambda_1,\ldots,\lambda_M$
satisfy {\bf L1} and {\bf L2}. If $\log(n \vee N) = o(n \wedge N)$,
then for $\nu=1,\ldots,M$,
\begin{equation}\label{eq:OPCA_risk_bound}
\sup_{\theta_\nu \in \mathbb{S}^{N-1}} \mathbb{E}L(\widehat
\theta_\nu, \theta_\nu) =
    \left[\frac{N-M}{nh(\lambda_\nu)} + \frac{1}{n}\sum_{\mu \neq \nu}
\frac{(\lambda_\mu+1)(\lambda_\nu+1)}{(\lambda_\nu -
  \lambda_\mu)^2}\right](1+o(1)).
\end{equation}}

\vskip.1in \noindent{\bf Remark :} It is possible to relax some of
the conditions stated in the theorem. On the other hand, with some
reasonable assumptions on the decay of the eigenvalues, it is also
possible to incorporate cases where $M$ is no longer a constant, but
increases with $n$. Then the issues would include, rates of growth
of $M$ and the rate of decay of eigenvalues that would result in the
OPCA estimator retaining consistency and the expression for its
asymptotic risk. These issues are not going to be addressed here.
However, it is important to note that, such questions have been
investigated - not necessarily for the Gaussian case - in the
context of spectral decomposition of $L^2$ stochastic processes by,
Hall and Horowitz (2004), Tony Cai and Hall (2005), Boente and
Fraiman (2000), Hall and Hosseini-Nasab (2006) among others.
However, these analyses do not deal with measurement errors. The
condition $\frac{N}{nh(\lambda_\nu)} \to 0$ is a necessary condition
for uniform convergence, as shown in Theorem 2.
It should be noted that, there are results, proved under slightly
different circumstances, that obtain the rates given by
(\ref{eq:OPCA_risk_bound}) as an upper bound on the rate of
convergence of OPCA estimators (Bai (2003), Cardot (2000), Kneip
(1994)). These analyses, while treating the problem under less
restrictive assumptions than Gaussianity (essentially, finite eighth
moment for the noise $Z_{ik}$), make the assumption that
$\frac{N^2}{n} \to 0$, when the $\lambda_\nu$'s are considered
fixed.

\section{Sparse model for eigenvectors}\label{sec:highdpca-sparsity}

In this section we discuss the concept of sparsity of the
eigenvectors and impose some restrictions on the space of
eigenvectors that lead to a sparse parametrization. This notion will
be used later from a decision-theoretic view point in order to
analyze the risk behavior of estimators of the eigenvectors. From
now on, $\theta$ will be used to denote the matrix
$[\theta_1,\ldots,\theta_M]$.

\subsection{$l^q$ constraint on the parameters}

The parameter space is taken to be a class of $M$-dimensional
positive semi-definite matrices satisfying the following criteria:
\begin{itemize}
\item $\lambda_1 > \ldots > \lambda_M$.
\item  For each $\nu=1,\ldots,M$, $\theta_\nu \in \Theta_\nu$
  for some $\Theta_\nu \subset \mathbb{S}^{N-1}$ that
  gives a sparse parametrization, in that most of the coefficients
  $\theta_{\nu k}$ are close to zero.
\item $\theta_1, \ldots, \theta_M$ are orthonormal.
\end{itemize}
One way to formalize the requirement of sparsity is to demand, as in
Johnstone and Lu (2004), that $\theta_\nu$ belongs to a weak-$l^q$
space $wl^q(C)$ where $C, q >0$. This space is defined as follows.
Suppose that the coordinates of a vector $\mathbf{x} \in
\mathbb{R}^N$ are $|x|_{(1)}, \ldots, |x|_{(N)}$, where $|x|_{(k)}$
is the $k$-th largest element, in absolute value. Then
\begin{equation}\label{eq:waek_l_q}
\mathbf{x} \in wl^q(C) ~~~~ \Leftrightarrow ~~~~ |x|_{(k)} \leq C
k^{-1/q}, ~~k=1,2,\ldots .
\end{equation}
In the Functional Data Analysis context, one can think of the
observations as the vectors of wavelet coefficients (when
transformed in an orthogonal wavelet basis of sufficient regularity)
of the observed functions. If the smoothness of a function $g$ is
measured by its membership in a Besov space $B_{q',r}^\alpha$, and
if the vector of its wavelet coefficients, when expanded in a
sufficiently regular wavelet basis, is denoted by $\mathbf{g}$, then
from Donoho (1993),
$$
g \in B_{q',r}^\alpha ~~~ \Longrightarrow ~~~  \mathbf{g} \in wl^q,
\qquad q = \frac{2}{2\alpha + 1}, ~~\mbox{if}~~\alpha > (1/q' -
1/2)_+.
$$
One may refer to Johnstone (2002) for more details. Treating this as
a motivation, instead of imposing a weak-$l^q$ constraint on the
parameter $\theta_\nu$, we rather impose an $l^q$ constraint. Note
that, for $C, q > 0$,
\begin{equation}\label{eq:l_q_constr}
\mathbf{x} \in \mathbb{R}^N \cap l^q(C)  ~~~~\Leftrightarrow ~~~~
\sum_{k=1}^N |x_k|^q \leq C^q.
\end{equation}
Since $l^q(C) \hookrightarrow wl^q(C)$, it is possible to derive
lower bounds on the minimax risk of estimators when the parameter
lies in a $wl^q$ space by restricting attention to an $l^q$ ball of
appropriate radius.

For $C > 0$, define $\Theta_q(C)$ by
\begin{equation}\label{eq:Theta_q_C}
\Theta_q(C) = \{ a \in \mathbb{S}^{N-1} : \sum_{k=1}^N |a_k|^q \leq
C^q\},
\end{equation}
where $\mathbb{S}^{N-1}$ is the unit sphere in $\mathbb{R}^N$
centered at 0.  One important fact is, if $0< q < 2$, for
$\Theta_q(C)$ to be nonempty, one needs $C \geq 1$, while for $q >
2$, the reverse inequality is necessary. Further, for $0 < q <2$, if
$C^q \geq N^{1-q/2}$, then the space $\Theta_q(C)$ reduces to
$\mathbb{S}^{N-1}$ because in this case, the vector
$(1/\sqrt{N},1/\sqrt{N},\ldots,1/\sqrt{N})$ is in the parameter
space. Also, the only vectors that belong in the space when $C=1$
are the poles, i.e. vectors of the form $(0,0,\ldots,0,\pm
1,0,\ldots,0)$, where the non-zero term appears in exactly one
coordinate. Define, for $q \in (0, 2)$,  $m_C$ to be an integer
$\geq 1$ that satisfies
\begin{equation}\label{eq:m_C_def}
m_C^{1-q/2} \leq C^q < (m_C+1)^{1-q/2} .
\end{equation}
Then $m_C$ is the largest dimension of a unit sphere, centered at 0,
that fits inside the parameter space $\Theta_q(C)$.

\subsection{Parameter space}\label{subsec:highdpca-parameter}

The parameter space for $\theta := [\theta_1:\ldots:\theta_M]$ is
denoted by
\begin{equation}\label{eq:Theta_q_M}
\Theta_q^M(C_1,\ldots,C_M) = \{ \theta \in \prod_{\nu=1}^M
\Theta_q(C_\nu) ~:~ \langle \theta_\nu, \theta_{\nu'} \rangle = 0,
~~\mbox{for}~~ \nu \neq \nu'\},
\end{equation}
where $\Theta_q(C)$ is defined through (\ref{eq:Theta_q_C}), and
$C_\nu \geq 1$ for all $\nu=1,\ldots,M$.

\vskip.1in\noindent{\bf Remark :} \label{rem:alternative_sparsity}
If $M > 1$, one can describe the sparsity of the eigenvectors in a
different way. Consider the sequence $\zeta := \zeta^N =
(\sqrt{\sum_{\nu=1}^M \lambda_\nu \theta_{\nu k}^2} :
k=1,2,\ldots,N)$. One may demand that the vector $\zeta$ be sparse
in an $l^q$ or weak-$l^q$ sense. This particular approach to
sparsity has some natural interpretability, since the quantity
$\zeta_k^2 = \sum_{\nu=1}^M \lambda_\nu \theta_{\nu k}^2$, where
$\zeta_k$ is the $k$-th coordinate of $\zeta$, is the variance of
the $k$-th coordinate of the ``signal'' part of the vector $X$.
There is a connection between this model and the model we intend to
study. If (\ref{eq:Theta_q_M}) holds, then $\zeta \in
l^q_N(\overline C_\lambda)$, where $\overline C_\lambda^q =
\sum_{\nu=1}^M \lambda_\nu^{q/2} C_\nu^q$. On the other hand, $l^q$
(weak-$l^q$) sparsity of $\zeta$ implies $l^q$ (weak-$l^q$) sparsity
of $\theta_\nu$ for all $\nu=1,\ldots,M$.

\subsection{Lower bound on the minimax risk}
\label{subsec:highdpca-lower_bound}

In this section a lower bound on the minimax risk of estimating
$\theta_\nu$ over the parameter space (\ref{eq:Theta_q_M}) is
derived when $0 < q < 2$, under the loss function defined through
(\ref{eq:loss_fn}). The result is stated under some simplifying
assumptions that make the asymptotic analysis more transparent.
Define
\begin{equation}\label{eq:g_lambda_tau}
g(\lambda,\tau) = \frac{(\lambda - \tau)^2}{(1+\lambda)(1+\tau)},
\qquad \lambda, \tau > 0.
\end{equation}
\begin{itemize}
\item[{\bf A1}] There exists a constant $C_0 > 0$ such that $C_0^q <
  C_\mu^q - 1$ for all $\mu=1,\ldots,M$, for all $N$.
\item[{\bf A2}] As $n, N \to \infty$, $n h(\lambda_\nu) \to \infty$.
\item[{\bf A3}] As $n, N \to \infty$, $n h(\lambda_\nu) = O(1)$.
\item[{\bf A4}] As $n, N \to \infty$, $n g(\lambda_\mu,\lambda_\nu)
\to \infty$ for all $\mu =1,\ldots,\nu-1,\nu+1,\ldots,M$.
\item[{\bf A5}] As $n, N \to \infty$, $n \max_{1\leq \mu \neq \nu \leq
    M} g(\lambda_\mu,\lambda_\nu) = O(1)$.
\end{itemize}
Conditions {\bf A4} and {\bf A5} are applicable only when $M > 1$.
In the statement of the following theorem, the infimum is taken over
all estimators $\widehat \theta_\nu$, estimating $\theta_\nu$,
satisfying $\parallel \widehat \theta_\nu \parallel = 1$.

\vskip.15in \noindent{\bf Theorem 2:}
\label{thm:theta_nu_minimax_lbd} \textit{Let $0 < q < 2$ and $1\leq
\nu \leq M$. Suppose that {\bf A1} holds.
\begin{itemize}
\item[(a)]
If {\bf A3} holds, then there exists $B_1 > 0$ such that
\begin{equation}\label{eq:theta_nu_minimax_lbd_a}
\liminf_{n \to \infty} \inf_{\widehat \theta_\nu} \sup_{\theta \in
\Theta_q^M(C_1,\ldots,C_M)} \mathbb{E} L(\widehat \theta_\nu,
\theta_\nu) \geq B_1 .
\end{equation}
\item[(b)] If {\bf A2} holds, then there exists $B_2 > 0$, $A_q > 0$,
  and $c_1 \in (0,1)$, such that
\begin{equation}\label{eq:theta_nu_minimax_lbd_b}
\liminf_{n \to \infty} \delta_n^{-1} \inf_{\widehat \theta_\nu}
\sup_{\theta \in \Theta_q^M(C_1,\ldots,C_M)} \mathbb{E} L(\widehat
\theta_\nu, \theta_\nu) \geq B_2,
\end{equation}
where $\delta_n$ is defined by
\begin{equation}\label{eq:theta_nu_minimax_lbd_b_delta_n}
\delta_n =
\begin{cases}
c_1 & ~~if~~ nh(\lambda_\nu) \leq \min\{c_1 (N-M), A_q
\overline{C}_\nu^q (n h(\lambda_\nu))^{q/2}\} \cr c_1
\frac{N-M}{nh(\lambda_\nu)} & ~~if~~ c_1 (N-M) \leq
\min\{nh(\lambda_\nu), A_q \overline{C}_\nu^q (n
h(\lambda_\nu))^{q/2}\} \cr A_q
\frac{\overline{C}_\nu^q}{(nh(\lambda_\nu))^{1-q/2}} & ~~if~~ A_q
\overline{C}_\nu^q (n h(\lambda_\nu))^{q/2} \leq
\min\{nh(\lambda_\nu), c_1 (N-M)\} \cr
\end{cases}
\end{equation}
and
\begin{equation}\label{eq:theta_nu_minimax_lbd_b_delta_n_alpha}
\delta_n = (c_2(\alpha))^{1-q/2} \frac{\overline{C}_\nu^q (\log
N)^{1-q/2}} {(nh(\lambda_\nu))^{1-q/2}}, ~~\mbox{if}~~A_{q,\alpha}
\overline{C}_\nu^q (\frac{n h(\lambda_\nu)}{\log N})^{q/2} \leq
\min\{\frac{nh(\lambda_\nu)}{\log N},K N^{1-\alpha}\},
\end{equation}
for some $K > 0$, $\alpha \in (0,1)$, $c_q(\alpha) \in (0,1)$ and
$A_{q,\alpha} > 0$. Here $\overline{C}_\nu^q := C_\nu^q - 1$. Also,
one can take $c_1 = \log(9/8)$, $A_q = (\frac{9 c_1 }{2})^{1-q/2}$,
$A_{q,\alpha} = (\alpha/2)^{1-q/2}$, $c_q(\alpha) =
(\alpha/9)^{1-q/2}$, $B_2 = \frac{1}{8}$ and $B_3 = (8e)^{-1}$.
\item[(c)] Suppose that $M > 1$. If
{\bf A4} holds, then there exists $B_3 >0$ such that
\begin{equation}\label{eq:theta_nu_minimax_lbd_c}
\liminf_{n \to \infty} \overline{\delta}_n^{-1} \inf_{\widehat
\theta_\nu} \sup_{\theta \in \Theta_q^M(C_1,\ldots,C_M)}\mathbb{E}
L(\widehat \theta_\nu, \theta_\nu) \geq B_3,
\end{equation}
where
\begin{equation}\label{eq:theta_nu_minimax_lbd_c_delta_n}
\overline{\delta}_n = \frac{1}{n} \max_{\mu
  \in\{1,\ldots,M\}\setminus\{\nu\}}
\frac{1}{g(\lambda_\mu,\lambda_\nu)}~.
\end{equation}
One can take $B_3 = \frac{1}{8e}$. However, if {\bf A5} holds, then
(\ref{eq:theta_nu_minimax_lbd_a}) is true.
\end{itemize}}

\vskip.1in \noindent{\bf Remark :} In the statement of Theorem 2,
there is much flexibility in terms of
what values the ``hyperparameters'' $C_1,\ldots, C_M$ and the
eigenvalues $\lambda_1,\ldots,\lambda_M$ can take. In particular,
they can vary with $N$, subject to the modest requirement that {\bf
A1} is satisfied. However, the constants appearing in equations
(\ref{eq:theta_nu_minimax_lbd_b}) and
(\ref{eq:theta_nu_minimax_lbd_c}) are not optimal.

\vskip.1in \noindent{\bf Remark :}\label{rem:minimax_lbd_size}
Another notable aspect is that, as the proof later shows, the rate
lower bounds in Part (b) are all of the form
$\frac{m}{nh(\lambda_\nu)}$, where $m$ is the ``effective'' number
of ``significant'' coordinates. This phrase becomes clear if one
notices further that, in the construction that leads to the lower
bound (see Section \ref{subsec:highdpca-minimax_lbd_b}), the vector
$\theta_\nu$ in a near-worst case scenario has overwhelming number
of coordinates of size $const. ~\frac{1}{\sqrt{nh(\lambda_\nu)}}$,
or, in the case (\ref{eq:theta_nu_minimax_lbd_b_delta_n_alpha}), of
size $const. ~\frac{\sqrt{\log N}}{\sqrt{nh(\lambda_\nu)}}$. Here
$m$ is of the same order as the number of these ``significant''
coordinates. This suggests that, an estimation strategy that is able
to extract coordinates of $\theta_\nu$ of the stated size, would
have the right rate of convergence, subject to possibly some
regularity conditions. The estimator described later (ASPCA) is
constructed by following this principle.

Part (a) and the second statement of Part (c) of Theorem 2
depict situations under which there is no estimator that is
asymptotically uniformly consistent over
$\Theta_q^M(C_1,\ldots,C_M)$. Moreover,  the first part of Part (b),
and Theorem 1
readily yield the following corollary.

\vskip.1in \noindent{\bf Corollary 1:}\label{cor:OPCA_minimax_rate}
\textit{If the conditions of Theorem 1
hold, and if {\bf A1} holds, together with the condition that
$$
\liminf_{n \to \infty} \frac{\overline{C}_\nu^q
  (nh(\lambda_\nu))^{q/2}}{N} > c_1^{q/2} A_q^{-1},
$$
then the usual PCA-based estimator of $\widehat \theta_\nu$, i.e.
the eigenvector corresponding to the $\nu$-th largest eigenvalue of
$\mathbf{S}$, has asymptotically the best rate of convergence.}

\vskip.1in\noindent{\bf Remark :} A closer look at the proof of
Theorem 1
reveals that the method of proof
explicitly made use of condition {\bf L1} to ensure that the
contribution of $\lambda_1,\ldots,\lambda_M$ to the residual term of
the second order expansion of $\widehat \theta_\nu$ is bounded.
However, the condition $n\max_{\mu \neq \nu}
g(\lambda_\mu,\lambda_\nu) \to \infty$ is certainly much weaker than
that. The method of proof pursued here fails to settle the question
as to whether this is sufficient to get the asymptotic rate
(\ref{eq:OPCA_risk_bound}). It is conjectured that this is the case.

\section{Estimation scheme}\label{sec:aspca-estimation}

This section outlines an estimation strategy for the eigenvectors
$\theta_\nu$, $\nu=1,\ldots,M$. Model (\ref{eq:basic}) is assumed
throughtout for observations $X_i$, $i=1,\ldots,n$. We propose
estimators is for the case when the noise variance $\sigma^2$ is
known. Therefore, without loss of generality, it can be taken to be
1. Henceforth, for simplicity of notations, it is also assumed that
$\xi = 0$. In practice, one may have to estimate $\sigma^2$ from
data. The median of the diagonal entries of the sample covariance
matrix $\mathbf{S} := \frac{1}{n}\mathbf{X}\mathbf{X}^T$ serves as a
reasonable (although slightly biased) estimator of $\sigma^2$, if
the true model is sparse. In the latter case, the data are rescaled
by multiplying each observation by $\widehat \sigma ^{-1}$, and the
resultant covariance matrix is called, with a slight abuse of
notation, $\mathbf{S}$. Note that, in this case, the estimates of
eigenvalues of $\Sigma$ are $\widehat\sigma^2$ times the
corresponding eigenvalues of $\mathbf{S}$.

\subsection{Sparse Principal Components Analysis
  (SPCA)}\label{subsec:aspca-estimation-spca}

In order to motivate the approach that is described in what follows,
consider first the SPCA estimation scheme studied by Johnstone and
Lu (2004). To that end, let $\mathbf{S} = \frac{1}{n}
\mathbf{X}\mathbf{X}^T$ denote the sample covariance matrix. Suppose
that the sample variances of coordinates (i.e., diagonal terms of
$\mathbf{S}$) are denoted by $\hat \sigma_1^2, \ldots,\hat
\sigma_N^2$.

\begin{itemize}
\item
Define $\widehat I_n$ to be the set of indices $k \in
\{1,\ldots,N\}$ such that $\hat \sigma_k^2 > \gamma_n$ for some
threshold $\gamma_n > 0$.
\item
Let $\mathbf{S}_{\widehat I_n,\widehat I_n}$ be the submatrix of
$\mathbf{S}$ corresponding to the coordinates $\widehat I_n$.
Perform an eigen-analysis of  $\mathbf{S}_{\widehat I_n, \widehat
I_n}$. Denote the eigenvectors by
$\mathbf{e}_1,\ldots,\mathbf{e}_{\min\{n,|\widehat I_n|\}}$.
\item
For $\nu=1,\ldots,M$, estimate $\theta_\nu$ by $\widetilde
{\mathbf{e}}_\nu$ where $\widetilde {\mathbf{e}}_\nu$, an $N \times
1$ vector, is obtained from $\mathbf{e}_\nu$ by augmenting zeros to
all the coordinates that are in $\{1,\ldots,N\} \setminus \widehat
I_n$.
\end{itemize}
Johnstone and Lu (2004) showed that, if one  chooses an appropriate
threshold $\gamma_n$, then the estimate of $\theta_\nu$ is
consistent under the weak-$l^q$ sparsity constraint on $\theta_\nu$.
However, Paul and Johnstone (2004) showed that even with the best
choice of $\gamma_n$, the rate of convergence of the risk of this
estimate is not optimal. Indeed, Paul and Johnstone (2004)
demonstrate an estimator which has a better rate of convergence in
the single component ($M =1$) situation.

\subsection{Augmented Sparse PCA
  (ASPCA)}\label{subsec:aspca-estimation-aspca}

We now propose the ASPCA estimation scheme. This scheme is a
refinement of the SPCA scheme of Johnstone and Lu (2004), and can be
viewed as a generalization of the estimation scheme proposed by Paul
and Johnstone (2004) in the single component ($M=1$) case.

The key idea behind this estimation scheme is that, in addition to
using the coordinates having large variance, if one also uses the
covariance structure appropriately, then under the assumption of a
sparse structure of the eigenvectors, one will be able to extract a
lot more information and thereby get more accurate estimate of the
eigenvalues and eigenvectors. Notice that SPCA only focuses on the
diagonal of the covariance matrix and therefore ignores the
covariance structure. This renders this scheme suboptimal from an
asymptotic minimax risk analysis point of view. To make this point
clearer, it is instructive to analyze the covariance matrix in the
$M=1$ case. In view of the second Remark
after the statement of Theorem 2
one expects to be able to recover coordinates $k$ for which
$|\theta_{1k}| \gg \frac{1}{\sqrt{nh(\lambda_1)}}$. However, the
best choice for $\gamma_n$ for SPCA is $\gamma \sqrt{\frac{\log
n}{n}}$, for some constant $\gamma > 0$, which is way too large. On
the other hand, suppose that one divides the coordinates into two
sets $A$ and $B$, where the former contains all those $k$ such that
$|\theta_k|$ is ``large'', and the latter contains smaller
coordinates. Partition the matrix $\Sigma$ as
$$
\Sigma = \begin{bmatrix} \Sigma_{AA} & \Sigma_{AB} \cr
                         \Sigma_{BA} & \Sigma_{BB} \cr
         \end{bmatrix}
$$
Here $\Sigma_{BA} = \lambda_1 \theta_{1,B} \theta_{1,A}^T$. Assume
that, there is a ``preliminary'' estimator of $\theta_1$, say
$\widetilde \theta_1$ such that, $\langle \widetilde \theta_{1,A},
\theta_{1,A}\rangle \to 1$ in probability as $n \to \infty$. Then
one can use this estimator as a ``filter'', in a way described
below, to recover the ``informative ones'' among the smaller
coordinates. This can be seen from the following relationship
$$
\Sigma_{BA} \widetilde \theta_{1,A} = \langle \widetilde
\theta_{1,A}, \theta_{1,A}\rangle \lambda_1 \theta_{1,B} \approx
\lambda_1 \theta_{1,B}.
$$
In this manner one can extract some information about those
coordinates of $\theta_1$ that are in set $B$. The algorithm
described below is a generalization of this idea. It has three
stages. First two stages will be referred to as ``coordinate
selection'' stages. The final stage consists of an eigen-analysis of
the submatrix of $\mathbf{S}$ corresponding to the selected
coordinates, followed by a hard thresholding of the estimated
eigenvectors.

Let $\gamma_i > 0$ for $i=1,2,3$ and $\kappa > 0$  be four constants
to be specified later. Define $\gamma_{1,n} = \gamma_1
\sqrt{\frac{\log(n \vee N)}{n}}$.

\begin{itemize}
\item[{$1^o$}]
Select coordinates $k$ such that $\widehat \sigma_{kk} :=
\mathbf{S}_{kk} > 1 + \gamma_{1,n}$. Denote the set of selected
coordinates by $\widehat I_{1,n}$.
\item[{$2^o$}]
Perform  spectral decomposition of $\mathbf{S}_{\widehat
I_{1,n},\widehat I_{1,n}}$. Denote the eigenvalues by $\widehat
\ell_1 > \ldots > \widehat \ell_{m_1}$ where $m_1 =
\min\{n,|\widehat I_{1,n}|\}$, and corresponding eigenvectors by
$\mathbf{e}_1,\ldots,\mathbf{e}_{m_1}$.
\item[{$3^o$}]
Estimate $M$ by $\widehat M$ defined in Section
\ref{subsec:aspca-M_hat}.  Estimate $\lambda_j$ by $\widetilde
\lambda_j = \widehat \ell_j-1$, $j=1,\ldots,\widehat M$.
\item[{$4^o$}]
Define $E = [\frac{1}{\sqrt{\widehat \ell_1}}\mathbf{e}_1 : \ldots :
\frac{1}{\sqrt{\widehat \ell_{\widehat
      M}}}\mathbf{e}_{\widehat M}]$.
Compute $\mathbf{Q} = \mathbf{S}_{\widehat I_{1,n}^c,\widehat
I_{1,n}} E$.
\item[{$5^o$}]
Denote the  diagonal of the matrix $\mathbf{Q}\mathbf{Q}^T$ by $T$.
Define   $\widehat I_{2,n}$ to be the set of coordinates $k \in
\{1,\ldots,N\} \setminus \widehat I_{1,n}$ such that $|T_k| >
\gamma_{2,n}^2$ where
$$
\gamma_{2,n} = \gamma_2 \left(\sqrt{\frac{\log(n \vee N)}{n} } +
\frac{1}{\kappa} \sqrt{\frac{\widehat M}{n}}\right).
$$
\item[{$6^o$}]
Take the union   $\widehat I_n := \widehat I_{1,n} \bigcup \widehat
I_{2,n}$. Perform  spectral decomposition of $\mathbf{S}_{\widehat
I_n,\widehat
  I_n}$ . Estimate
$\theta_\nu$ by  augmenting the $\nu$-th eigenvector, with zeros in
the coordinates $\{1,\ldots,N\} \setminus \widehat I_n$, for
$\nu=1,\ldots, \widehat M$. Call this vector $\widehat \theta_\nu$.
\item[{$7^o$}]
Perform a coordinatewise ``hard'' thresholding of $\widehat
\theta_\nu$ at threshold
$$
\gamma_{3,n} := \gamma_3 \sqrt{\frac{\log (n \vee N)}{nh(\widetilde
\lambda_\nu)}} ~,
$$
and then normalize the thresholded vectors to get the final estimate
$\overline \theta_\nu$.
\end{itemize}

\vskip.1in\noindent{\bf Remark :} \label{rem:gamma_choice} The
scheme is specified except for the ``tuning parameters''
$\gamma_1$,$\gamma_2$,$\gamma_3$ and $\kappa$. The choice of
$\gamma_i$'s is discussed in the context of deriving upper bounds on
the risk of the estimator. It will be shown that, it suffices to
take $\gamma_1 = 4$, $\kappa = 2+\epsilon$ for a small $\epsilon >
0$, and $\gamma_2 = \sqrt{\frac{3}{2}} \kappa$. An analysis of the
thresholding scheme is not done here, but in practice $\gamma_3 = 3$
works well enough, and some calculations suggest that $\gamma_3 = 2$
suffices asymptotically.

\subsection{Estimation of $M$}\label{subsec:aspca-M_hat}

Let $\overline{\gamma}_1, \gamma_1' > 0$ be such that
$\overline{\gamma}_1  >  \gamma_1'$. Define
\begin{eqnarray}
\widehat{\overline{I}}_{1,n} &=& \{ k : \mathbf{S}_{kk} >
1+\overline{\gamma}_{1,n} \} ~~\mbox{where}~~
\overline{\gamma}_{1,n} = \overline{\gamma}_1
\sqrt{\frac{\log (n \vee N)}{n}}~, \label{eq:gamma_1_bar_n} \\
\widehat I_{1,n}' &=& \{ k : \mathbf{S}_{kk} > 1+\gamma_{1,n}'\}
~~\mbox{where}~~ \gamma_{1,n}' = \gamma_1' \sqrt{\frac{\log (n \vee
    N)}{n}}~.     \label{eq:gamma_1_prime_n}
\end{eqnarray}
Define
\begin{equation}\label{eq:M_alpha_n}
\alpha_n = 2\sqrt{\frac{|\widehat I_{1,n}'|}{n}} + \frac{|\widehat
  I_{1,n}'|}{n}
+ 6\left(\frac{|\widehat
  I_{1,n}'|}{n} \vee 1\right) \sqrt{\frac{\log (n \vee |\widehat
  I_{1,n}'|)}{n \vee |\widehat I_{1,n}'|}} ~.
\end{equation}
Let $\widehat{\overline \ell}_1 > \ldots > \widehat{\overline
  \ell}_{\overline{m}_1}$, where $\overline{m}_1 =
\min\{n,|\widehat{\overline I}_{1,n}|\}$, be the nonzero eigenvalues
of $\mathbf{S}_{\widehat{\overline I}_{1,n},\widehat{\overline
    I}_{1,n}}$. Define $\widehat M$ by
\begin{equation}\label{eq:M_hat_def}
\widehat M = \max \{1\leq k \leq \overline{m}_1 : \widehat{\overline
  \ell}_k  > 1+\alpha_n\}.
\end{equation}
The choice of $\gamma_1'$ and $\overline{\gamma}_1$ is discussed in
Section \ref{subsec:aspca-consistency_M_hat}.

\vskip.1in\noindent{\bf Remark :}\label{rem:non_sparse} Sparsity of
the eigenvectors is an implicit assumption for ASPCA scheme.
However, in practice, and specifically with only moderately large
samples, it is not always the case that ASPCA is able to select the
significant coordinates. More importantly, the scheme produces a
\textit{bona fide} estimator only when $\widehat I_{1,n}$ is
non-empty. If this is not the case, then one may use the $\nu$-th
eigenvector of $\mathbf{S}$ as the estimator of $\theta_\nu$.
However, determination of $M$ in this situation is a difficult
issue, and without recourse to additional information, one may set
$\widehat M = 0$.

\section{Rates of convergence}\label{sec:aspca-rates}

In this section we describe the asymptotic risk of ASPCA estimators
under some regularity conditions. The risk is analyzed under the
loss function (\ref{eq:loss_fn}), and it is assumed that condition
{\bf BA} of Section \ref{sec:model} holds. Further, the parameter
space for $\theta = [\theta_1:\ldots:\theta_M]$, over which the risk
is maximized, is taken to be $\Theta_q^M(C_1,\ldots,C_M)$ defined
through (\ref{eq:Theta_q_M}) in Section
\ref{subsec:highdpca-parameter},  where $0 < q <2$ and $C_1, \ldots,
C_M > 1$.

\subsection{Sufficient conditions for
  convergence}\label{subsec:aspca-rates-condition}

The following conditions are imposed on the ``hyperparameters'' of
the parameter space $\Theta_q(C_1,\ldots,C_M)$. Suppose that
$\rho_1,\ldots,\rho_M$ are as in {\bf C1} given below. Define
\begin{equation}\label{eq:rho_q_C}
\rho_q(C) := \sum_{\nu=1}^M \rho_\nu^{q/2} C_\nu^q.
\end{equation}
Observe that, since $C_\nu \geq 1$ for all $\nu=1,\ldots,M$,
$\rho_q(C) \geq \sum_{\nu=1}^M \rho_\nu^{q/2} \geq 1$.
\begin{itemize}
\item[{\bf C1}]
$\lambda_1,\ldots,\lambda_M$ are such that, as $n \to \infty$,
$\frac{\lambda_\nu}{\lambda_1} \to \rho_\nu$ where $1 \equiv \rho_1
> \rho_2
> \ldots > \rho_M$.
\item[{\bf C2}]
$\log N \asymp \log n$ and $\frac{(\log n)^2}{n\lambda_1^2} \to 0$
as $n \to \infty$.
\item[{\bf C3}]
$\frac{\rho_q(C) (\log N)^{1/2-q/4}}{\lambda_1^{1-q/2}
  n^{1/2-q/4}} \to 0$ as $n \to \infty$.
\end{itemize}
We discuss briefly the importance of these conditions. {\bf C1} is a
repetition of {\bf L1}. {\bf C2} is a convenient and very mild
technical assumption that should hold in most practical situations.
Second part of ${\bf C2}$ is non-trivial only when $\lambda_1 \to 0$
as $n \to \infty$. {\bf C3} requires some  explanation. It will
become increasingly clear that,  in order to get a uniformly
consistent estimate of the eigenvectors from the preliminary SPCA
step, one needs {\bf C3} to hold. Indeed, the sequence described in
{\bf C3} has the same asymptotic order as a common upper bound for
the rate of convergence of the supremum risk of the SPCA estimators
of all the $\theta_\nu$'s. So, the implication is that if {\bf C3}
holds then the SPCA scheme of Johnstone and Lu (2004) gives
consistent estimates.

\vskip.1in\noindent{\bf Remark :}\label{rem:ASPCA_rate} Note that,
$\frac{1}{nh(\lambda)} \leq \frac{1+c}{n\lambda^2}$ if $\lambda \in
(0,c)$ and $\frac{1}{nh(\lambda)} \leq \frac{1}{\eta(c)
  n\lambda}$ if
$\lambda \geq c$, for any $c > 0$. Since $\rho_q(C) \geq 1$, {\bf
C3} guarantees that
\begin{equation}\label{eq:ASPCA_rate}
\frac{\rho_q(C) (\log N)^{1-q/2}}{(n
  h(\lambda_1))^{1-q/2}} = o(1), \qquad \mbox{as}~ n \to \infty.
\end{equation}
In fact, if $\liminf_{n\to \infty}  \lambda_1 \geq c > 0$, then the
upper bound in (\ref{eq:ASPCA_rate}) can be replaced by
$o((\frac{\log N}{n})^{1/2-q/4})$. It will be shown that this is a
common (and near-optimal) upper bound  on the rate of convergence of
the ASPCA estimate of $\theta_\nu$'s. If one compares this with the
lower bound given by Theorem 2,
it is conjectured that (\ref{eq:ASPCA_rate}) should also be a
sufficient condition for establishing that the lower bound defined
through (\ref{eq:theta_nu_minimax_lbd_b_delta_n_alpha}) is also the
upper bound on the minimax risk, at the level of rates. However,
since our method depends on finding a preliminary consistent
estimator of the eigenvectors (in our case SPCA), the somewhat
stronger condition {\bf C3} becomes necessary to establish rates of
convergence of the ASPCA estimator.

\subsection{Statement of the result}\label{subsec:aspca-result}

Now we state the main result of this section. The asymptotic
analysis of risk is conducted only for the estimator $\widehat
\theta_\nu$ for eigenvector $\theta_\nu$, and not for the
thresholding estimator $\widetilde \theta_\nu$. Derivation of the
results for $\widetilde \theta_\nu$ requires additional technical
work, but can be carried out. It can be shown that in certain
circumstances the latter has a slightly better asymptotic risk
property. In practice, the thresholding estimator seems to work
better when the eigenvalues are well-separated. The following
theorem describes the asymptotic behavior of the risk of the ASPCA
estimator $\widehat \theta_\nu$ under the loss function $L$ defined
through (\ref{eq:loss_fn}). $g(\cdot, \cdot)$ is defined by
(\ref{eq:g_lambda_tau}).

\vskip.15in\noindent{\bf Theorem 3:}\label{thm:ASPCA_minimax_bound}
\textit{Assume that {\bf BA} and conditions {\bf C1}-{\bf C3} hold.
Then, there are constants $K:=K(q,\gamma_1,\gamma_2,\kappa)$ and
$K':= K'(q,M,\gamma_1,\gamma_2, \kappa)$ such that, as $n \to
\infty$, for all $\nu=1,\ldots,M$,
\begin{eqnarray}\label{eq:ASPCA_minimax_bound}
&&\hskip-.3in \sup_{\theta \in \Theta_q^M(C_1,\ldots,C_M)}
\mathbb{E}L(\widehat
\theta_{\nu},\theta_{\nu}) \nonumber\\
&\leq& \hskip-.1in\left[K (C_{\nu}^q + K'\rho_\nu^{-q}
\frac{\rho_q(C)}{\log (n\vee N)}) \left(\frac{\log (n \vee
    N)}{nh(\lambda_{\nu})}\right)^{1-q/2} + \sum_{\mu\neq
  \nu}^M \frac{1}{ng(\lambda_\mu,\lambda_\nu)}\right](1+o(1))
\end{eqnarray}}

\vskip.1in\noindent{\bf Remark :}\label{rem:ASPCA_minimaxity} The
expression in the upper bound is somewhat cumbersome, but the
significance of each of the terms in (\ref{eq:ASPCA_minimax_bound})
will become clear in the course of the proof. However, notice that,
if the parameters $C_1,\ldots,C_M$ of the space
$\Theta_q(M)(C_1,\ldots,C_M)$ are such that,
\begin{equation}\label{eq:C_nu_codition}
\exists~0 < \underline{C} < \overline{C} < \infty, ~~~\mbox{such
  that}~~~ \underline{C} \leq
\frac{\max_{1\leq\mu\leq M} C_\mu}{\min_{1\leq \mu \leq M} C_\mu}
\leq \overline{C}, ~~~\mbox{for all}~~ n,
\end{equation}
then, Theorem 3
and Theorem 2
together imply that, under conditions {\bf BA}, {\bf C1}-{\bf C3},
{\bf A1} and the condition on the hyperparameters given by
(\ref{eq:theta_nu_minimax_lbd_b_delta_n_alpha}), the ASPCA estimator
$\widehat \theta_\nu$ has the optimal rate of convergence. The
condition (\ref{eq:C_nu_codition}) is satisfied in particular if
$C_1,\ldots,C_M$ are all bounded above.

\vskip.1in It is important to emphasize that
(\ref{eq:ASPCA_minimax_bound}) is an asymptotic result in the
following sense. It is possible to give finite a sample bound on
$\sup_{\theta \in \Theta_q^M(C_1,\ldots,C_M)} \mathbb{E}L(\widehat
\theta_\nu,\theta_\nu)$. However, this upper bound involves many
additional terms whose total contribution is smaller than a
prescribed $\epsilon > 0$ only when $n \geq n_\epsilon$, say, where
$n_\epsilon$ depends on the hyperparameters, apart from $\epsilon$.

\vskip.1in\noindent{\bf Remark :}\label{rem:rate_comparison_OPCA} It
is instructive to compare the asymptotic supremum risk of ASPCA with
that of OPCA (or usual PCA based) estimator of $\theta_\nu$. A
closer inspection of the proof reveals that, if for all sufficiently
large $n$,
\begin{equation*}
N \leq K'' (\frac{\rho_q(C)}{\log (n\vee N)}) \left(\frac{\log (n
\vee
    N)}{nh(\lambda_{\nu})}\right)^{-q/2},
\end{equation*}
then for some constant $K'' > 0$, under {\bf BA} and {\bf C1}-{\bf
C3}, one can replace the upper bound in
(\ref{eq:ASPCA_minimax_bound}) by
$$
\overline{K} \left[\frac{N\log (n \vee N)}{nh(\lambda_\nu)} +
\sum_{\mu\neq \nu}
\frac{1}{ng(\lambda_\mu,\lambda_\nu)}\right](1+o(1)),
$$
for some constant $\overline{K}$. This rate is greater than that of
OPCA estimator by a factor of at most $\log (n \vee N)$. However,
observe that, the bound on the risk of OPCA estimator holds under
weaker conditions. In particular,  Theorem 1
does not assume any particular structure for the eigenvectors.

\section{Proof of Theorem 2}
\label{sec:highdpca_proof_minimax_lbd}

The proof requires a closer look at the geometry of the parameter
space, in order to obtain good finite dimensional subproblems that
can then be used as inputs to the general machinery, to come up with
the final expressions.

\subsection{Risk bounding
  strategy}\label{subsec:highdpca_lower_bound_strategy}

A key tool for our proof the lower bound on the minimax risk is
Fano's lemma. Thus,  it is necessary to derive a general expression
for the Kullback-Leibler discrepancy between the probability
distributions described by two separate parameter values.

\vskip.15in\noindent{\bf Proposition 1:}\label{prop:multi_KL_div}
\textit{Let $\theta^{(j)} = [\theta_1^{(j)} : \ldots :
\theta_M^{(j)}]$, $j=1,2$ be two parameters. Let $\Sigma_{(j)}$
denote the matrix given by (\ref{eq:Sigma_basic}) with $\theta =
\theta^{(j)}$ (and $\sigma = 1$). Let $P_j$ denote the joint
probability distribution of $n$ i.i.d. observations from
$N(0,\Sigma_{(j)})$. Then the Kullback-Leibler discrepancy of $P_2$
from $P_1$, to be denoted by $K_{1,2} :=
K(\theta^{(1)},\theta^{(2)})$, is given by
\begin{equation}\label{eq:multi_KL_div}
K_{1,2} \stackrel{def}{=} K(\theta^{(1)},\theta^{(2)}) = n\left[
\frac{1}{2}\sum_{\nu=1}^M \eta(\lambda_\nu) \lambda_\nu -
\frac{1}{2} \sum_{\nu =1}^M \sum_{\nu'=1}^M \eta(\lambda_\nu)
\lambda_{\nu'}  |\langle \theta_{\nu'}^{(1)},
\theta_\nu^{(2)}\rangle|^2\right],
\end{equation}
where
\begin{equation}\label{eq:eta_lambda}
\eta(\lambda) = \frac{\lambda}{1+\lambda}, \qquad \lambda > 0.
\end{equation}}

\subsection{Use of Fano's lemma}

We outline the general approach pursued in the rest of this section.
The idea is to bound the supremum of the risk on the entire
parameter space by the maximum risk over a finite subset of it, and
then to use some variant of Fano's lemma to provide a lower bound
for the latter quantity.

Thus, the goal is to find an appropriate \textit{finite} subset
${\cal F}_0$ of $\Theta_q^M(C_1,\ldots,C_M)$, such that the
following properties hold.
\begin{itemize}
\item[(1)]{If $\theta^{(1)}, \theta^{(2)} \in {\cal F}_0$,
then $L(\theta_\nu^{(1)},\theta_\nu^{(2)}) \geq 4\delta$, for some
    $\delta > 0$ (to be chosen). This property will be referred to as
  ``$4\delta$-distinguishability in $\theta_\nu$''.}
\item[(2)]{The element $\theta \in {\cal F}_0$ is a unique
    representative of the equivalence class $[\theta]$, where
    $[\theta]$ is defined to be the class of $N \times M$
    matrices whose $\nu$-th column is either $\theta_\nu$ or $-\theta_\nu$.}
\item[(3)]{Subject to (1), the quantity $\sup_{i \neq j :
      ~\theta^{(i)},
\theta^{(j)} \in {\cal F}_0}
    K(\theta^{(i)},\theta^{(j)}) + K(\theta^{(j)},\theta^{(i)})$ is as
    small as possible.}
\end{itemize}
Given any estimator $\widehat \theta$ of $\theta$, based on data
$\mathbf{X}_n = (X_1,\ldots,X_n)$, define a new estimator
$\phi(\mathbf{X}_n)$ (an $N \times M$ matrix) as $\phi(\mathbf{X}_n)
= \theta^*$ if $\theta^* = \arg \min_{\theta \in {\cal F}_0}
L(\theta_\nu,\widehat \theta_\nu)$, where $\widehat \theta_\nu$ is
the $\nu$-th column of $\widehat \theta$ (i.e., estimate of
$\theta_\nu$). Then, by Chebyshev's inequality,
\begin{eqnarray}\label{eq:risk_lb_M}
\sup_{\theta \in \Theta_q^M(C_1,\ldots,C_M)} \mathbb{E}_\theta
L(\theta_\nu,\widehat \theta_\nu) &\geq&  \delta \sup_{\theta \in
\Theta_q^M(C_1,\ldots,C_M)} \mathbb{P}_\theta (
L(\theta_\nu,\widehat \theta_\nu) \geq \delta)  \nonumber\\
&\geq& \delta \sup_{\theta \in {\cal F}_0} \mathbb{P}_\theta (
L(\theta_\nu,\widehat \theta_\nu) \geq \delta) \nonumber\\
&\geq& \delta \sup_{\theta \in {\cal F}_0} \mathbb{P}_\theta (
[\phi(\mathbf{X}_n)] \neq [\theta]).
\end{eqnarray}
The last inequality is because, if $L(\theta_\nu^{(j)},\widehat
\theta_\nu) < \delta$ for any $\theta^{(j)} \in {\cal F}_0$, then by
the ``$4\delta$-distinguishability in $\theta_\nu$'' (property (1)
above), it follows that $[\phi_\nu(\mathbf{X}_n)] =
[\theta_\nu^{(j)}]$,  and hence $[\phi(\mathbf{X}_n)] =
[\theta^{(j)}]$.

Two versions of Fano's lemma are found to be useful in this context.
The following version, due to Birg\'{e} (2001), of a result of Yang
and Barron (1999) (p.1570-71), is most suitable when ${\cal F}_0$
can be chosen to be large.

\vskip.15in\noindent{\bf Lemma 1:} \label{lemma:Fano_lb_large}
\textit{Let $\{P_\theta : \theta \in \Theta\}$ be a family  of
probability distributions on a common measurable space, where
$\Theta$ is an arbitrary parameter space. Suppose that a loss
function for the estimation problem is given by $L'(\theta,\theta')
= \mathbf{1}_{\theta \neq \theta'}$. Define the minimax risk over
$\Theta$ by
$$
p_{max} = \inf_T  \sup_{\theta \in \Theta} \mathbb{P}_\theta(T \neq
\theta), = \inf_T  \sup_{\theta \in \Theta} \mathbb{E} L'(\theta,T),
$$
where $T$ denotes an arbitrary estimator of $\theta$ with values in
$\Theta$. Then for any finite subset ${\cal F}$ of $\Theta$, with
elements $\theta_1,\ldots,\theta_J$ where $J = |{\cal F}|$,
\begin{equation}\label{eq:Fano}
p_{max} \geq 1 -  \inf_{Q} ~\frac{J^{-1} \sum_{i=1}^J K(P_i, Q) +
\log
  2}{\log J}
\end{equation}
where $P_i =  \mathbb{P}_{\theta_i}$, and $Q$ is an arbitrary
probability distribution, and $K(P_i,Q)$ is the Kullback-Leibler
divergence of $Q$ from $P_i$.}

\vskip.15in To use Lemma 1
choose $P_i$ to
be $P_{\Sigma_{(i)}} \equiv P_{\theta^{(i)}} := N^{\otimes
  n}(0,\Sigma_{(i)})$, where  $\Sigma_{(i)}$ is the matrix
$\sum_{\nu=1}^M \lambda_\nu \theta_\nu^{(i)} {\theta_\nu^{(i)}}^T +
I$, and $\theta^{(i)} \in {\cal F}_0$ $i=1,\ldots, |{\cal F}_0|$,
are the distinct values of parameter $\theta$ that constitute the
set ${\cal F}_0$. Then set $Q_0 = P_{\theta^{(0)}}$, for some
appropriately chosen $\theta^{(0)} \in \Theta_q^M(C_1,\ldots,C_M)$
such that the following condition is satisfied.
\begin{equation}\label{eq:nicely_spaced}
ave_{1 \leq i \leq |{\cal F}_0|} K(\theta^{(i)},\theta^{(0)}) \asymp
\sup_{1 \leq i \leq |{\cal F}_0|} K(\theta^{(i)},\theta^{(0)}),
\end{equation}
where the notation ``$\asymp$'' means that the both sides are are
within constant multiples of each other. Then it follows from
(\ref{eq:risk_lb_M}) and Lemma 1
that,
\begin{equation}\label{eq:risk_lbd}
\delta^{-1} \sup_{\theta \in \Theta_q^M(C_1,\ldots,C_M)}
\mathbb{E}_\theta L(\theta_\nu,\widehat \theta_\nu) \geq 1 -
\frac{ave_{1\leq i \leq |{\cal F}_0|} K(\theta^{(i)},\theta^{(0)}) +
\log 2} {\log |{\cal F}_0|}~.
\end{equation}
To complete the picture it is desirable that
\begin{equation}\label{eq:Fano_lb_approx}
\frac{ave_{1\leq i \leq |{\cal F}_0|} K(\theta^{(i)},\theta^{(0)}) +
\log 2} {\log |{\cal F}_0|} ~\approx~ c,
\end{equation}
where $c$ is a number between 0 and 1.

A different version of Fano's lemma, due to Birg\'{e} (2001), is
needed when ${\cal F}_0$ consists of only two elements
$\theta^{(1)}$ and $\theta^{(2)}$, so that the classification
problem reduces to a test of hypothesis of $P_1$ against $P_2$.

\vskip.15in\noindent{\bf Lemma 2:}\label{lemma:Fano_lb_small}
\textit{Let $\alpha_T$ and $\beta_T$ denote respectively the Type I
and Type II errors associated with an arbitrary test $T$ between the
two simple hypotheses $P_1$ and $P_2$. Define, $\pi_{mis} = \inf_T
(\alpha_T + \beta_T)$, where the infimum is taken over all test
procedures.
\begin{equation}\label{eq:risk_lb_small}
K(P_1,P_2) \geq - \log[\pi_{mis}(2-\pi_{mis})].
\end{equation}}

\subsection{Geometry of the parameter space}\label{subsec:highdpca-geometry}

We view the space $\Theta_q(C)$, for $0 < q <2$, as the
$N$-dimensional unit sphere centered at the origin, from which some
parts have been chopped off, symmetrically in each coordinate, such
that there is some portion left at each pole (i.e., a point of the
form $(0,\ldots,0,\pm 1,0,\ldots,0)$, where the non-zero term
appears only once).
In this connection, we define an object that is central to the proof
of Theorem \ref{thm:theta_nu_minimax_lbd}.

\vskip.1in\noindent{\bf Definition :} \label{defn:m_r_polar_sphere}
\textit{Let $0 < r < 1$ and $N
> m \geq 1$. An $(N,m,r)$ {\bf polar sphere} at pole $k_0$, on set
$J = \{j_1,\ldots,j_m\}$, where $1\leq k_0 \leq N$ and $j_l \in
\{1,\ldots,N\}\setminus \{k_0\}$ for $l=1,\ldots,m$, is a subset of
~$\mathbb{S}^{N-1}$ given by
\begin{equation}\label{eq:N_m_r_polar_sphere}
{\cal S}(N,m,r,k_0,J) := \{ \mathbf{x} \in \mathbb{S}^{N-1} :
x_{k_0} = \sqrt{1-r^2}, \sum_{l=1}^m x_{j_l}^2 = r^2 \}.
\end{equation}}
So, an $(N,m,r)$ polar sphere is centered at the point
$(0,\ldots,0,\sqrt{1-r^2},0,\ldots,0)$, (which is not in
$\mathbb{S}^{N-1}$), has radius $r$, and has dimension $m$. Note
that, the largest sphere of any given dimension $m$, such that $C^q
< m^{1-q/2}$ (equivalently, $m > m_C$, where $m_C$ is defined
through (\ref{eq:m_C_def})), that can be inscribed inside
$\Theta_q(C)$ is an $(N,m,r)$ polar sphere. The radius $r$ of such a
polar sphere, given $C^q < m^{1-q/2}$ (or $m > m_C$), to be denoted
by $r_m(C)$, satisfies
\begin{equation}\label{eq:polar_radius}
\{1-(r_m(C))^2\}^{q/2} + m^{1-q/2} \{r_m(C)\}^q = C^q.
\end{equation}
Of course, if $C^q \geq m^{1-q/2}$ (or $m_C \geq m$ ) then as a
convention, $r_m(C) = 1$. Condition (\ref{eq:polar_radius}) ensures
that all the points lying on an $(N,m,r)$ polar sphere such that $r
\in (0,r_m(C))$, are inside $\Theta_q(C)$.

\subsection{A common recipe for Part (a) and Part
  (b)}\label{subsec:highdpca-recipe}

In the proof of Part (a) and Part (b) of the theorem, there is a
common theme in the construction of ${\cal F}_0$. Let
$\mathbf{e}_\mu$ denote the $N$-vector whose $\mu$-th coordinate is
1 and rest are all zero. In either case, if $\{\theta^{(j)},
j=1,\ldots,|{\cal F}_0|\}$ is an enumeration of the elements of
${\cal F}_0$, then the following are true.
\begin{itemize}
\item
[(F1)] There is an $N \times M$ matrix $\theta^{(0)}$, such that
$\theta_\nu^{(0)} = \mathbf{e}_\nu$.
\item
[(F2)] $\theta_\mu^{(j)} = \mathbf{e}_\mu$ for $\mu =1,
\ldots,\nu-1,\nu+1,\ldots,M$, for all $j=0,1,\ldots,|{\cal F}_0|$.
\item
[(F3)] $\theta_\nu^{(j)} \in {\cal S}(N,m,r,\nu,J)$ for some $m$,
$r$ and $J$.  $m$ and $r$ are fixed for all $1\leq j \leq |{\cal
F}_0|$, but $J$ may be different for different $j$, depending on the
situation.
\end{itemize}
The $\theta^{(0)}$ in (F1) is the same $\theta^{(0)}$ appearing in
(\ref{eq:risk_lbd}). Also, (\ref{eq:multi_KL_div}) simplifies to
\begin{equation}\label{eq:KL_div}
K(\theta^{(j)},\theta^{(0)}) = \frac{1}{2} nh(\lambda_\nu) (1 -
(\langle \theta_\nu^{(j)},\theta_\nu^{(0)}\rangle)^2) = \frac{1}{2}
nh(\lambda_\nu) r^2, ~~j=1,\ldots,|{\cal F}_0|.
\end{equation}
Moreover, in either case, the points $\theta^{(j)}$ are so chosen
that
\begin{equation}\label{eq:r_square_distinguish}
L(\theta_\nu^{(j)},\theta_\nu^{(k)}) \geq r^2, ~~~\mbox{for all}~~~
1 \leq j\neq k \leq |{\cal F}_0|.
\end{equation}
In other words, the set ${\cal F}_0$ is $r^2$ distinguishable in
$\theta_\nu$.

\subsection{Proof of Part (a)}\label{subsec:highdpca-minimax_lbd_a}

Construct ${\cal F}_0$ satisfying (F1)-(F3), with
$$
\theta_\nu^{(j)} = \sqrt{1-r^2} \mathbf{e}_\nu + r \mathbf{e}_j,
\qquad j=M+1, \ldots, N,
$$
where $r \in (0,1)$ is such that $(1-r^2)^{q/2} + r^q \leq C_\nu^q$.
Thus, $|{\cal F}_0| = N - M$. Verify that
(\ref{eq:r_square_distinguish}) holds, in fact the lower bound is
$2r^2$, with an equality. Therefore, (\ref{eq:risk_lbd}) applies,
with $\delta = \frac{r^2}{2}$. Since $nh(\lambda_\nu)$ is bounded
above, and $\log (N-M) \to \infty$ as $n \to \infty$,
(\ref{eq:theta_nu_minimax_lbd_a}) follows from (\ref{eq:KL_div}).

\subsection{Connection to ``Sphere
  packing''}\label{subsec:highdpca-sphere_packing}

Our proof of Part (b) of Theorem 2
depends crucially on the following construction due to Zong (1999).

Let $m$ be a large positive integer, and  $m_0 =
\left[\frac{2m}{9}\right]$ (the largest integer $\leq
\frac{2m}{9}$). Define $Y_m^*$ as the maximal set of points of the
form ${\bf z}=(z_1,\ldots,z_m)$ in $\mathbb{S}^{m-1}$ such that the
following is true.
\begin{equation}\label{eq:separation}
\sqrt{m_0} z_i \in \{-1,0,1\} ~\forall~ i, ~~ \sum_{i=1}^m |z_i| =
\sqrt{m_0} ~~\mbox{and, for} ~~{\bf z}, {\bf z}' \in Y_m^*, ~~
\parallel {\bf z} - {\bf z}'\parallel \geq 1.
\end{equation}
For any $m \geq 1$, the maximal number of points lying on
$\mathbb{S}^{m-1}$ such that any two points are at distance at least
1, is exactly same as the \textit{kissing number} of an $m$-sphere.
It is known that this number is $\leq 3^m$ and $\geq
(9/8)^{m(1+o(1))}$. Zong (1999) uses the construction described
above to derive the lower bound, by showing that $|Y_m^*| \geq
(9/8)^{m(1+o(1))}$ for $m$ large.

\subsection{Proof of  Part (b)}\label{subsec:highdpca-minimax_lbd_b}

Structures of ${\cal F}_0$ for the three cases in
(\ref{eq:theta_nu_minimax_lbd_b_delta_n}) are similar. Set $m \leq
(N-M)$, large. Set $c_1 = \log(9/8)$, $A_q = (9c_1/2)^{1-q/2}$.
Choose $r \approx \sqrt{\delta_n}$, and define the set ${\cal F}_0$
satisfying (F1)-(F3) and the following construction.

Set $|{\cal F}_0| = |Y_m^*|$, where $Y_m^*$ is the set defined in
Section \ref{subsec:highdpca-sphere_packing}. Set,
\begin{equation}\label{eq:theta_nu_j_F_0}
\theta_\nu^{(j)} = \sqrt{1-r^2} \mathbf{e}_\nu + r \sum_{l=1}^m
z_l^{(j)} \mathbf{e}_{l+M}, ~~~j=1,\ldots,|{\cal F}_0|,
\end{equation}
where $\mathbf{z}^{(j)}= (z_1^{(j)},\ldots,z_m^{(j)})$, $j\geq 1$,
is an enumeration of the elements of $Y_m^*$. Observe that, for all
$j \geq 1$,
\begin{equation}\label{eq:theta_nu_j_on_polar}
\theta_\nu^{(j)} \in {\cal  S}(N,m,r,\nu,\{M+1,\ldots,M+m\}) \bigcap
{\cal S}(N,m_0,r,\nu,supp(\mathbf{z}^{(j)})),
\end{equation}
where $supp(\mathbf{z}^{(j)})$ is the set of nonzero coordinates of
$\mathbf{z}^{(j)}$. Therefore, (\ref{eq:r_square_distinguish}) and
(\ref{eq:KL_div}) hold for all $j \geq 1$.

\subsubsection{Case : $nh(\lambda_\nu) \leq \min\{c_1 (N-M), A_q
\overline{C}_\nu^q (n h(\lambda_\nu))^{q/2}\}$}

Take $m = [nh(\lambda_\nu)]$ and $r^2 = c_1$. Observe that, for all
$j \geq 1$,
$$
\parallel \theta_\nu^{(j)} \parallel_q^q = (1-r^2)^{q/2} + m_0^{1-q/2}
r^q \leq 1 + (2/9)^{1-q/2} (nh(\lambda_\nu))^{1-q/2} c_1^{q/2} \leq
1+c_1\overline{C}_\nu^q < C_\nu^q.
$$
Thus, ${\cal F}_0 \subset \Theta_q^M(C_1,\ldots,C_M)$. Further,
since $nh(\lambda_\nu) \to \infty$, $\log |{\cal F}_0| \geq c_1
nh(\lambda_\nu)(1+o(1))$. Since (\ref{eq:r_square_distinguish}) and
(\ref{eq:KL_div}) hold, with $\delta = \frac{r^2}{4}$, from
(\ref{eq:risk_lbd}) the result follows, because
$$
\limsup_{n\to \infty} \frac{ave_{1\leq j \leq |{\cal F}_0|}
K(\theta^{(j)},\theta^{(0)}) +
  \log 2}{\log |{\cal F}_0|} \leq \limsup_{n \to \infty}
\frac{\frac{1}{2} c_1 nh(\lambda_\nu) + \log 2}{ c_1
nh(\lambda_\nu)} =\frac{1}{2}~.
$$

\subsubsection{Case : $c_1 (N-M)
\leq \min\{nh(\lambda_\nu), A_q \overline{C}_\nu^q (n
h(\lambda_\nu))^{q/2}\}$}

Take $m = N-M$ and $r^2 = \frac{c_1(N-M)}{nh(\lambda_\nu)}$. Then,
for all $j \geq 1$,
$$
\parallel \theta_\nu^{(j)} \parallel_q^q \leq 1
+ (2/9)^{1-q/2} (N-M) c_1^{q/2} (nh(\lambda_\nu))^{-q/2} \leq
1+\overline{C}_\nu^q = C_\nu^q.
$$
The result follows by arguments similar to those used for the case
$nh(\lambda_\nu) \leq \min\{c_1 (N-M), A_q \overline{C}_\nu^q (n
h(\lambda_\nu))^{q/2}\}$.

\subsubsection{Case : $A_q \overline{C}_\nu^q (n
h(\lambda_\nu))^{q/2} \leq \min\{nh(\lambda_\nu), c_1 (N-M) \}$}

Take $m = [c_1^{-q/2} (9/2)^{1-q/2} \overline{C}_\nu^q (n
h(\lambda_\nu))^{q/2}]$ and $r^2 = c_1 \frac{m}{nh(\lambda_\nu)}$.
Again, verify that $m \to \infty$ as $n \to \infty$ (by {\bf A1}),
and for $j \geq 1$,
$$
\parallel \theta_\nu^{(j)} \parallel_q^q \leq 1
+ (2/9)^{1-q/2} m^{1-q/2} c_1^{q/2}
(\frac{m}{nh(\lambda_\nu)})^{q/2} \leq 1+\overline{C}_\nu^q =
C_\nu^q,
$$
and the result follows by familiar arguments.

\subsubsection{Proof of (\ref{eq:theta_nu_minimax_lbd_b_delta_n_alpha})}

The construction in all three previous cases assumes that the set of
non-zero coordinates is held fixed (in our case
$\{M+1,\ldots,M+m\}$) for every fixed $m$. However, it is possible
to get a bigger set ${\cal F}_0$ satisfying the requirements, if
this condition is relaxed.

Suppose that $A_{q,\alpha} = (\alpha/2)^{1-q/2}$, and the condition
in (\ref{eq:theta_nu_minimax_lbd_b_delta_n_alpha}) holds for some
$\alpha \in (0,1)$. Set $m =
[(\alpha/9)^{-q/2}(9/2)^{1-q/2}\overline{C}_\nu^q (n
h(\lambda_\nu))^{q/2}(\log N)^{-q/2}]$ and $r^2 = (\alpha/9)
\frac{m}{nh(\lambda_\nu)}$. Take $c_q(\alpha) = (\alpha/9)^{1-q/2}$.
Observe that $m \to \infty$ as $n \to \infty$, $m = O(N^{1-\alpha})$
and $r \in (0, 1)$.  Set $\theta^{(0)} = [\mathbf{e}_1:\ldots
:\mathbf{e}_M]$. For every set $\pi \subset \{M+1,\ldots,N\}$ of
size $m$, construct ${\cal F}_\pi$ satisfying (F1)-(F3) such that,
\begin{equation}\label{eq:theta_nu_j_F_pi}
\theta_\nu^{(j)} = \sqrt{1-r^2} \mathbf{e}_\nu + r \sum_{l \in \pi}
z_l^{(j)} \mathbf{e}_l, ~~~j=1,\ldots,|Y_m^*|.
\end{equation}
As before,  ${\cal F}_\pi \subset \Theta_q^M(C_1,\ldots,C_M)$, for
all $\pi$, so that (\ref{eq:KL_div}) and
(\ref{eq:r_square_distinguish}) are satisfied. Let ${\cal P}$ to be
a collection of such sets $\pi$ such that, for any two sets $\pi$
and $\pi'$ in ${\cal P}$, the set $\pi \cap \pi'$ has cardinality at
most $\frac{m_0}{2}$. This ensures that
$$
\mbox{for} ~~ \mathbf{y}, \mathbf{y}' \in \bigcup_{\pi \in \cal P}
{\cal F}_\pi, \quad L(\mathbf{y},\mathbf{y}') \geq r^2.
$$
This also ensures that the sets ${\cal F}_\pi$ are disjoint  for
$\pi \neq \pi'$, since each $\theta_\nu^{(j)}$ for $\theta^{(j)} \in
{\cal F}_0$ is nonzero in exactly $m_0+1$ coordinates. Define ${\cal
F}_0 = \bigcup_{\pi \in \cal P} {\cal F}_\pi$. Then
\begin{equation}\label{eq:big_F_0_lower}
|{\cal F}_0| = |\bigcup_{\pi \in \cal P} {\cal F}_\pi| = |{\cal P}|
~|Y_m^*| \geq |{\cal P}| (9/8)^{m(1+o(1))}.
\end{equation}
By Lemma 7,
stated in Section \ref{subsec:highdpca-counting}, there is a
collection ${\cal P}$ such that $|{\cal P}|$ is at least $\exp( [N
{\cal E}(m/9N) - 2m{\cal E}(1/9)](1+o(1)))$, where ${\cal E}(x)$ is
the Shannon entropy function :
$$
{\cal E}(x) = - x \log(x) - (1-x) \log(1-x), ~~ 0 < x < 1.
$$
Since ${\cal E}(x) \sim -x\log x$ when $x \to 0+$, it follows from
(\ref{eq:big_F_0_lower}) that,
$$
\frac{\log |{\cal F}_0|}{m} \geq  [ \frac{1}{9}(\log N - \log m) - 2
{\cal
  E}(1/9) + \log 9
+ \log(9/8)](1+o(1)) \geq \frac{\alpha}{9} \log N (1+o(1)),
$$
since $m = O(N^{1-\alpha})$. Finally, observe that
$$
\limsup_{n\to \infty} \frac{ave_{\theta^{(j)} \in  |{\cal F}_0|}
K(\theta^{(j)},\theta^{(0)}) +
  \log 2}{\log |{\cal F}_0|} \leq \limsup_{n \to \infty}
\frac{\frac{1}{2} (\alpha/9) m \log N}{ (\alpha/9) m \log N}
=\frac{1}{2}
$$
and use (\ref{eq:risk_lbd}) to finish argument.

\subsection{Proof of Part (c)}\label{subsec:highdpca-minimax_lbd_c}

Consider first the proof of (\ref{eq:theta_nu_minimax_lbd_c}). Fix a
$\mu \in \{1,\ldots,M\} \setminus \{\nu\}$. Define $\theta^{(1)}$
and $\theta^{(2)}$ as follows. Set $r^2 =
\frac{2}{ng(\lambda_1,\lambda_2)}$ (assume w.l.o.g. that $r < 1
\wedge C_0$). Take $\theta_{\mu'}^{(j)} = \mathbf{e}_{\mu'}$,
$j=1,2$ for all $\mu' \neq \mu,\nu$. Define
\begin{equation}\label{eq:optimal_theta_mu}
\theta_\nu^{(1)} = \mathbf{e}_\nu, ~~~\theta_\nu^{(2)} =
\sqrt{1-r^2} \mathbf{e}_\nu + r\mathbf{e}_\mu, ~~~ \theta_\mu^{(1)}
= \mathbf{e}_\mu, ~~ \theta_\mu^{(2)} = -r\mathbf{e}_\nu +
\sqrt{1-r^2} \mathbf{e}_\mu.
\end{equation}
Observe that $\theta_\mu^{(j)}\perp \theta_\nu^{(j)}$, $j=1,2$,
$\langle\theta_\nu^{(1)},\theta_\nu^{(2)}\rangle = \sqrt{1-r^2} =
\langle\theta_\mu^{(1)},\theta_\mu^{(2)}\rangle$ and
$\langle\theta_\mu^{(1)},\theta_\nu^{(2)}\rangle = r = -
\langle\theta_\nu^{(1)},\theta_\mu^{(2)}\rangle$. Also, by {\bf A1},
$\theta^{(j)} \in \Theta_q^M(C_1,\ldots,C_M)$, for $j=1,2$.

Let $P_j = N^{\otimes n}(0,\Sigma_{(j)})$. Then
\begin{eqnarray}\label{eq:2_KL_div_sym}
K(P_1,P_2) + K(P_2,P_1) &=& n [  h(\lambda_\mu) (1- |\langle
\theta_\mu^{(1)}, \theta_\mu^{(2)}\rangle|^2) + h(\lambda_\nu) (1 -
|\langle\theta_\nu^{(1)},\theta_\nu^{(2)}\rangle|^2)
\nonumber\\
&& ~~~ - \frac{1}{2}(\lambda_{\mu}\eta(\lambda_\nu) +
\lambda_{\nu}\eta(\lambda_\mu)) \{|\langle \theta_\mu^{(1)},
\theta_\nu^{(2)}\rangle|^2 +|\langle \theta_\nu^{(1)},
\theta_\mu^{(2)}\rangle|^2 \} ]
\nonumber\\
&=& n[(h(\lambda_\mu) + h(\lambda_\nu))r^2 -
\frac{1}{2}(\lambda_{\mu}\eta(\lambda_\nu) +
\lambda_{\nu}\eta(\lambda_\mu))r^2 ]\nonumber\\
&=& n g(\lambda_\mu,\lambda_\nu) r^2.
\end{eqnarray}
Apply Lemma 2
for testing $P_1$ against
$P_2$. Define $p_{mis} = \inf_T (\alpha_T \vee \beta_T)$ and observe
that $p_{mis} \leq \pi_{mis} \leq 2p_{mis}$.  Since the lower bound
in (\ref{eq:risk_lb_small}) is symmetric w.r.t. $\pi_{mis}$, and
$\pi_{mis}$ is symmetric w.r.t. $P_1$ and $P_2$, it follows that
$$
ng(\lambda_\mu,\lambda_\nu)r^2 = K(P_1,P_2) + K(P_2,P_1) \geq -2\log
(\pi_{mis}(2-\pi_{mis})).
$$
This implies that
$$
e^{-\frac{n}{2}g(\lambda_\mu,\lambda_\nu)r^2} \leq
\pi_{mis}(2-\pi_{mis}) \leq 2\pi_{mis} \leq 4p_{mis}
$$
Since, $L(\theta^{(1)},\theta^{(2)}) = 2(1-\sqrt{1-r^2}) \geq r^2$,
and $r^2 = \frac{2}{ng(\lambda_\mu,\lambda_\nu)}$, use
(\ref{eq:risk_lb_M}) with ${\cal F}_0 =\{\theta^{(1)},
\theta^{(2)}\}$ and $\delta = r^2$ to get,
$$
\sup_{\theta \in \Theta_q^M(\theta_1,\ldots,\theta_M)}
\mathbb{E}_\theta L(\theta_\nu,\widehat \theta_\nu) \geq
\frac{1}{8e} \frac{1}{ng(\lambda_\mu,\lambda_\nu)}.
$$
Now, let $\mu$ vary over all the indices
$1,\ldots,\nu-1,\nu+1,\ldots,M$ and the result follows.

In the situation where $\overline{\delta}_n \not\to 0$, as $n \to
\infty$, simply take $\mu$ ($\neq \nu$) to be the index for which
$g(\lambda_\mu,\lambda_\nu)$ is minimum. Then apply the same
procedure as in above with $r \in (0,C_0)$ fixed.

\section{Proof of Theorem 1}
\label{subsec:highdpca-OPCA_risk_bound}

We require two main tools in the proof of Theorem 1
- one (Lemma 5) is concerned with the deviations of the extreme
eigenvalues of a Wishart$(N,n)$ matrix and the other (Lemma 6)
relates to the change in the eigen-structure of a symmetric matrix
caused by a small, additive perturbation. Sections
\ref{subsec:highdpca-eigen_deviation_bound} and
\ref{subsec:highdpca-perturb_eigen} are devoted to them. The
importance of Lemma 6
is that, in order to bound the risk of an estimator of $\theta_\nu$
one only needs to compute the expectation of squared norm of a
quantity that is linear in $\mathbf{S}$ (or a submatrix of this, in
case of ASPCA estimator). The second bound in
(\ref{eq:eigenvec_error}) then ensures that the remainder is
necessarily of smaller order of magnitude. This fact is used
explicitly in deriving (\ref{eq:L_theta_nu_theta_hat_nu_bound}).

\vskip.15in\noindent{\bf Remark :}\label{rem:H_nu} In view of Lemma
6,
$H_\nu(\Sigma)$ becomes a key quantity in the analysis of the risk
of any estimator of $\theta_\nu$. Observe that,
\begin{equation}\label{eq:H_nu}
H_\nu := H_\nu(\Sigma) = \sum_{1\leq \nu' \neq \nu \leq M}
\frac{1}{\lambda_{\nu'}
  -\lambda_\nu}\theta_{\nu'}\theta_{\nu'}^T - \frac{1}{\lambda_\nu}
(I - \sum_{\nu'=1}^M \theta_{\nu'}\theta_{\nu'}^T), \qquad
\nu=1,\ldots,M.
\end{equation}

\vskip.1in\noindent Expand matrix $\mathbf{S}$ as follows.
\begin{eqnarray}\label{eq:S_M_expand}
\mathbf{S} &=& \sum_{\mu=1}^M \frac{\parallel v_\mu \parallel^2}{n}
\lambda_\mu \theta_\mu \theta_\mu^T + \sum_{\mu=1}^M
\sqrt{\lambda_\mu} \left(\theta_\mu (\frac{1}{n} \mathbf{Z}
v_\mu)^T+ \frac{1}{n} \mathbf{Z}
  v_\mu \theta_\mu^T\right)  \nonumber\\
&& ~~~~~~ + \sum_{\mu \neq \mu'} \frac{\langle v_\mu,
v_{\mu'}\rangle}{n} \sqrt{\lambda_\mu \lambda_{\mu'}} \theta_\mu
\theta_{\mu'}^T + \frac{1}{n} \mathbf{Z}\mathbf{Z}^T .
\end{eqnarray}
In order to use Lemma 6,
an expression for $H_\nu \mathbf{S} \theta_\nu$ is needed. Use the
fact that $H_\nu \theta_\nu = 0$ and $\theta_\nu^T \theta_\mu =
\delta_{\mu \nu}$ (Kronecker's symbol), to conclude that
\begin{eqnarray}\label{eq:H_nu_S_theta_nu}
H_\nu \mathbf{S}\theta_\nu &=& \sum_{\mu\neq \nu}
\left(\sqrt{\lambda_\mu} \frac{1}{n}\langle \mathbf{Z} v_\mu,
\theta_\nu \rangle + \sqrt{\lambda_\mu \lambda_\nu} \frac{1}{n}
\langle v_\mu, v_\nu \rangle \right) H_\nu \theta_\mu \nonumber\\
&& ~~~~~ + \sqrt{\lambda_\nu} H_\nu \frac{1}{n} \mathbf{Z} v_\nu +
H_\nu \frac{1}{n} \mathbf{Z}\mathbf{Z}^T \theta_\nu.
\end{eqnarray}
Further, from (\ref{eq:H_nu}) it follows that, $H_\nu \theta_\mu =
\frac{1}{\lambda_\mu - \lambda_\nu} \theta_\mu$, if $\mu \neq \nu$.
Also,
\begin{equation}\label{eq:H_nu_Z_v_nu}
H_\nu \mathbf{Z} v_\nu = - \frac{1}{\lambda_\nu} (I -\sum_{\mu=1}^M
\theta_{\mu}\theta_{\mu}^T) \mathbf{Z} v_\nu + \sum_{\mu \neq \nu}
\frac{1}{\lambda_\mu - \lambda_\nu} \langle \mathbf{Z} v_\nu,
\theta_\mu \rangle \theta_\mu ,
\end{equation}
and
\begin{equation}\label{eq:H_nu_Z_Z_theta_nu}
H_\nu \mathbf{Z}\mathbf{Z}^T \theta_\nu = -\frac{1}{\lambda_\nu}(I
-\sum_{\mu=1}^M \theta_{\mu}\theta_{\mu}^T) \mathbf{Z}\mathbf{Z}^T
\theta_\nu +\sum_{\mu \neq \nu} \frac{1}{\lambda_\mu - \lambda_\nu}
\langle \mathbf{Z}^T \theta_\mu, \mathbf{Z}^T\theta_\nu \rangle
\theta_\mu .
\end{equation}
From (\ref{eq:H_nu_S_theta_nu}), (\ref{eq:H_nu_Z_v_nu}) and
(\ref{eq:H_nu_Z_Z_theta_nu}), it follows that
\begin{eqnarray}\label{eq:H_nu_S_theta_nu_expand}
H_\nu \mathbf{S}\theta_\nu &=&  \sum_{\mu\neq \nu}
\frac{1}{\lambda_\mu
  - \lambda_\nu}\left(\sqrt{\lambda_\mu}
\frac{1}{n}\langle \mathbf{Z} v_\mu, \theta_\nu \rangle +
\sqrt{\lambda_\nu} \frac{1}{n}\langle \mathbf{Z} v_\nu, \theta_\mu
\rangle\right)\theta_\mu \nonumber\\
&&  ~~~ + \sum_{\mu\neq \nu} \frac{1}{\lambda_\mu
  - \lambda_\nu} \left(\sqrt{\lambda_\mu \lambda_\nu} \frac{1}{n}
\langle v_\mu, v_\nu \rangle + \frac{1}{n}\langle \mathbf{Z}^T
\theta_\mu, \mathbf{Z}^T \theta_\nu
\rangle \right) \theta_\mu \nonumber\\
&&  ~~~ - \frac{1}{n\lambda_\nu}(I -\sum_{\mu=1}^M
\theta_{\mu}\theta_{\mu}^T) \mathbf{Z}\mathbf{Z}^T \theta_\nu -
\frac{1}{n\sqrt{\lambda_\nu}} (I -\sum_{\mu=1}^M
\theta_{\mu}\theta_{\mu}^T) \mathbf{Z} v_\nu .
\end{eqnarray}
Let $\Gamma$ be an $N \times (N-M)$ matrix such that $\Gamma^T
\Gamma = I$, and $\Gamma \Gamma^T = (I -\sum_{\mu=1}^M
\theta_{\mu}\theta_{\mu}^T)$. Then, $\Gamma \theta_\mu = 0$ for all
$\mu=1,\ldots,M$.

\vskip.1in\noindent\label{rem:Z_v_indep} A crucial fact here is
that, since $v_\mu$ has i.i.d. $N(0,1)$ entries, and is independent
of $\mathbf{Z}$, for any $D \in \mathbb{R}^{m\times n}$,
$D\mathbf{Z} \frac{v_\mu}{\parallel v_\mu
  \parallel}$ has a $N_m(0,DD^T)$ distribution, and is independent of
$v_\mu$. Furthermore, since $\theta_\mu$ are orthonormal, and
$\Gamma \theta_\mu = 0$ for all $\mu$, it follows that
$\mathbf{Z}^T\theta_{\mu}$ has a $N_n(0,I)$ distribution;
$\{\mathbf{Z}^T\theta_{\mu}\}_{\mu=1}^M$ are mutually independent
and are independent of $\Gamma \mathbf{Z}$.

\vskip.1in

Next, we compute some expectations that will lead to the final
expression for $\mathbb{E}\parallel H_\nu \mathbf{S}
\theta_\nu\parallel^2$.
\begin{eqnarray}\label{eq:H_nu_S_theta_nu_expect_2}
&&\mathbb{E}\left(\sqrt{\lambda_\mu} \frac{1}{n}\langle \mathbf{Z}
v_\mu, \theta_\nu \rangle + \sqrt{\lambda_\nu} \frac{1}{n}\langle
\mathbf{Z} v_\nu, \theta_\mu
\rangle\right)^2 \nonumber\\
&=& \frac{1}{n^2}\left[\lambda_\mu \mathbb{E}(\langle \mathbf{Z}
  v_\mu, \theta_\nu \rangle)^2 + \lambda_\nu \mathbb{E}
(\langle \mathbf{Z} v_\nu, \theta_\mu\rangle)^2 + 2\sqrt{\lambda_\mu
  \lambda_\nu}\mathbb{E}(\langle \mathbf{Z}
  v_\mu, \theta_\nu \rangle \langle \mathbf{Z} v_\nu,
  \theta_\mu\rangle)\right]^2 \nonumber\\
&=& \frac{\lambda_\mu + \lambda_\nu}{n}~,
\end{eqnarray}
since the cross product term vanishes, which can be verified by a
simple conditioning argument. By similar calculations,
\begin{equation}\label{eq:H_nu_S_theta_nu_expect_3}
\mathbb{E} \left(\sqrt{\lambda_\mu \lambda_\nu} \frac{1}{n} \langle
v_\mu, v_\nu \rangle + \frac{1}{n}\langle \mathbf{Z}^T \theta_\mu,
\mathbf{Z}^T \theta_\nu \rangle \right)^2 =
\frac{\lambda_\nu\lambda_\mu + 1}{n}~,
\end{equation}
and
\begin{equation}\label{eq:H_nu_S_theta_nu_expect_4}
\mathbb{E}\left(\sqrt{\lambda_\mu} \frac{1}{n}\langle \mathbf{Z}
v_\mu, \theta_\nu \rangle + \sqrt{\lambda_\nu} \frac{1}{n}\langle
\mathbf{Z} v_\nu, \theta_\mu \rangle\right)\left(\sqrt{\lambda_\mu
\lambda_\nu} \frac{1}{n} \langle v_\mu, v_\nu \rangle +
\frac{1}{n}\langle \mathbf{Z}^T \theta_\mu, \mathbf{Z}^T \theta_\nu
\rangle \right) = 0 .
\end{equation}
Since $trace(\Gamma \Gamma^T) = N-M$, from the remark made above, it
follows that,
\begin{eqnarray}
\mathbb{E}\parallel (I -\sum_{\mu=1}^M \theta_{\mu}\theta_{\mu}^T)
\mathbf{Z}\mathbf{Z}^T \theta_\nu \parallel^2 &=&
\mathbb{E}[(\theta_\nu^T \mathbf{Z})\mathbf{Z}^T \Gamma\Gamma^T
\mathbf{Z}
(\mathbf{Z}^T \theta_\nu)] =  n(N -M),\label{eq:H_nu_S_theta_nu_expect_5}\\
\mathbb{E}\parallel (I -\sum_{\mu=1}^M
\theta_{\mu}\theta_{\mu}^T)\mathbf{Z}v_\nu
\parallel^2 &=& \mathbb{E}\parallel v_\nu \parallel^2 \mathbb{E}\parallel
\Gamma^T \mathbf{Z}\frac{v_\nu}{\parallel v_\nu\parallel}
\parallel^2 = n (N-M), \label{eq:H_nu_S_theta_nu_expect_6}
\end{eqnarray}
and
\begin{equation} \label{eq:H_nu_S_theta_nu_expect_7}
\mathbb{E}\langle (I -\sum_{\mu=1}^M \theta_{\mu}\theta_{\mu}^T)
\mathbf{Z}\mathbf{Z}^T \theta_\nu, (I -\sum_{\mu=1}^M
\theta_{\mu}\theta_{\mu}^T)\mathbf{Z}v_\nu \rangle = \mathbb{E}
[v_\nu^T \mathbf{Z}^T \Gamma \Gamma^T \mathbf{Z} (\mathbf{Z}^T
\theta_\nu)] = 0.
\end{equation}
Use (\ref{eq:H_nu_S_theta_nu_expand}), and equations
(\ref{eq:H_nu_S_theta_nu_expect_2}) -
(\ref{eq:H_nu_S_theta_nu_expect_7}), together with the
orthonormality of $\theta_\mu$'s and the fact that $\Gamma\theta_\mu
= 0$ for all $\mu$ to conclude that,
\begin{equation} \label{eq:H_nu_S_theta_nu_expect}
\mathbb{E}\parallel H_\nu \mathbf{S}\theta_\nu \parallel^2 =
\frac{N-M}{nh(\lambda_\nu)} + \frac{1}{n} \sum_{\mu \neq \nu}
\frac{(1+\lambda_\mu)(1+\lambda_\nu)}{(\lambda_\mu - \lambda_\nu)^2}
~.
\end{equation}
The next step in the argument is to show that, $\max_{0\leq \mu \leq
M} (\lambda_\mu - \lambda_{\mu+1})^{-1} \parallel \mathbf{S} -
\Sigma
\parallel$ is small with a very high probability. Here, by convention,
$\lambda_0 = \infty$ and $\lambda_{M+1} = 0$. From
(\ref{eq:S_M_expand}),
\begin{eqnarray}\label{eq:S_norm_bound}
\parallel \mathbf{S} - \Sigma \parallel &\leq& \sum_{\mu=1}^M \lambda_\mu
|\frac{\parallel v_\mu \parallel^2}{n} - 1| + 2\sum_{\mu=1}^M
\sqrt{\lambda_\mu} \frac{1}{n} \parallel \mathbf{Z}v_\mu \parallel \nonumber\\
&& ~~~~~ + \sum_{\mu \neq \mu'} \sqrt{\lambda_\mu \lambda_{\mu'}}
|\frac{\langle v_\mu,
  v_{\mu'}\rangle}{n}| + \parallel \frac{1}{n} \mathbf{Z}\mathbf{Z}^T
- I\parallel.
\end{eqnarray}
Define, for any $c > 0$, $D_{1,n}(c)$ to be the set
\begin{eqnarray}\label{eq:D_1_def}
D_{1,n}(c) &=& \bigcap_{\mu=1}^M \{|\frac{\parallel v_\mu
\parallel^2}{n}-1| \leq 2c \sqrt{\frac{\log (n \vee N)}{n}}\}
\nonumber\\
&& \bigcap~\bigcap_{\mu=1}^M\{ \frac{\parallel \mathbf{Z} v_\mu
  \parallel}{n} \leq \left(1+
2c \sqrt{\frac{\log (n \vee N)}{n \wedge N}}\right)
\sqrt{\frac{N}{n}} \}\nonumber\\
&& \bigcap ~\bigcap_{1\leq \mu < \mu'\leq M} \{|\frac{\langle v_\mu,
v_{\mu'}\rangle}{n}| \leq c \sqrt{\frac{\log (n \vee N)}{n}}\} .
\end{eqnarray}
Use Lemmas 14 and 15
to prove that,
\begin{equation}\label{eq:D_1_prob_bound}
1 - \mathbb{P}(D_{1,n}(c)) \leq 3M (n \vee N)^{-c^2} + M(M-1) (n
\vee N)^{-\frac{3}{2}c^2 + O(\log(n \vee N)/n)} .
\end{equation}
Define $D_{2,n}(c)$ as
\begin{equation}\label{eq:D_2_def}
D_{2,n}(c) = \{\parallel \frac{1}{n} \mathbf{Z}\mathbf{Z}^T -
I\parallel \leq  2\sqrt{\frac{N}{n}} + \frac{N}{n} + c t_n\},
\end{equation}
with $t_n$ as in Lemma 5.
From
(\ref{eq:S_norm_bound}), (\ref{eq:D_1_prob_bound}) and
(\ref{eq:eigen_deviation_bound}), it follows that for $n \geq n_c$,
\begin{eqnarray}\label{eq:S_norm_prob_bound}
&& \mathbb{P}(\parallel \mathbf{S} - \Sigma \parallel >
\epsilon_{n,N}(c,\lambda)) ~\leq~ 1-\mathbb{P}(D_{1,n}(c) \cap
D_{2,n}(c))
\nonumber \\
&\leq& (3M+2) (n \vee N)^{-c^2} + M(M-1) (n \vee N)^{-\frac{3}{2}c^2
+ O(\log(n \vee N)/n)},
\end{eqnarray}
where
\begin{eqnarray}\label{eq:epsilon_lambda_n_N}
\epsilon_{n,N}(c,\lambda) &=& 2c (\sum_{\mu=1}^M \lambda_\mu)
\sqrt{\frac{\log (n \vee N)}{n}} + 2(\sum_{\mu=1}^M
\sqrt{\lambda_\mu}) \left(1+ 2c \sqrt{\frac{\log (n \vee N)}{n
\wedge N}}\right)
\sqrt{\frac{N}{n}} \nonumber\\
&& ~~~+ c (\sum_{1\leq \mu \neq \mu'\leq M} \sqrt{\lambda_\mu
\lambda_{\mu'}} ) \sqrt{\frac{\log (n \vee N)}{n}} +
2\sqrt{\frac{N}{n}} + \frac{N}{n} + c t_n .
\end{eqnarray}
Define
\begin{equation}\label{delta_lambda_n_N}
\delta_{n,N,\nu} = \max \{ (\lambda_\nu - \lambda_{\nu+1})^{-1},
(\lambda_{\nu-1} - \lambda_\nu)^{-1}\}
\epsilon_{n,N}(\sqrt{2},\lambda),
\end{equation}
and observe that $\delta_{n,N,\nu} \to 0$ as $n \to \infty$ under
{\bf L1} and {\bf L2}.

To complete the proof of (\ref{eq:OPCA_risk_bound}), write
\begin{equation}\label{eq:OPCA_theta_nu_expand}
\widehat \theta_\nu - sign(\theta_\nu^T \widehat \theta_\nu)
\theta_\nu= - H_\nu \mathbf{S} \theta_\nu + R_\nu.
\end{equation}
Since $\delta_{n,N,\nu} \to 0$,  by (\ref{eq:eigenvec_error}),
(\ref{eq:eigen_Delta_r}), (\ref{eq:eigen_Delta_bar_r}) and
(\ref{eq:S_norm_prob_bound}), and the fact that $\Delta_r \leq
\overline{\Delta}_r$, for sufficiently large $n$, on
$D_{1,n}(\sqrt{2}) \cap D_{2,n}(\sqrt{2})$,
\begin{equation}\label{eq:L_theta_nu_theta_hat_nu_bound}
\parallel H_\nu \mathbf{S} \theta_\nu \parallel^2 (1-\delta_{n,N,\nu}')^2
\leq L(\theta_\nu, \widehat \theta_\nu)  \leq \parallel H_\nu
\mathbf{S} \theta_\nu \parallel^2 (1+\delta_{n,N,\nu}')^2,
\end{equation}
where
\begin{equation}
\delta_{n,N,\nu}' = \frac{\delta_{n,N,\nu}}
{(1-2\delta_{n,N,\nu}(1+2\delta_{n,N,\nu}))^2}
[1+2(1+\delta_{n,N,\nu})(1-2\delta_{n,N,\nu}(1+2\delta_{n,N,\nu}))],
\end{equation}
and $\delta_{n,N,\nu}' \to 0$ as $n \to \infty$. Since
$L(\theta_\nu, \widehat \theta_\nu) \leq 2$,
(\ref{eq:S_norm_prob_bound}),
(\ref{eq:L_theta_nu_theta_hat_nu_bound}) and
(\ref{eq:H_nu_S_theta_nu_expect}) together imply
(\ref{eq:OPCA_risk_bound}).

\section{Proof of Theorem 3}

In some respect the proof of Theorem 3
bears resemblance to the proof of Theorem 4 in Johnstone and Lu
(2004). The basic idea in both these cases is to first provide a
``bracketing relation''. This means that, if $\widehat I_n$ denotes
the set of selected coordinates, and $\underline{I}_n$ and
$\overline{I}_n$ are two \textit{non-random} sets with suitable
properties, then an inequality of the form
$\mathbb{P}(\underline{I}_n \subset \widehat I_n \subset
\overline{I}_n)  \geq 1 - b_n$ holds, where $b_n$ converges to zero
at least polynomially in $n$. Once this relationship is established,
one can utilize it to study the eigen-structure of the submatrix
$\mathbf{S}_{\widehat I_n, \widehat I_n}$ of $\mathbf{S}$. The
advantage of this is that the bracketing relation ensures that the
quantities involved in the perturbation terms for the eigenvectors
and eigenvalues can be controlled, except possibly on a set of
probability at most $b_n$.

The proof of Theorem 3
follows this principle. However, there are several technical aspects
in both the steps that require much computation. The first step,
namely, establishing a bracketing relation for $\widehat I_n$, is
done in Sections \ref{sec:aspca-bracketing} -
\ref{subsec:aspca-second_stage}. The second step follows more or
less the approach taken in the proof of Theorem 1,
in that, on a set of high probability, an upper bound on $L(\widehat
\theta_\nu,\theta_\nu)$ is established that is of the form
$\parallel H_\nu \mathbf{S}_2 \theta_\nu\parallel^2(1+\delta_n)$,
where $H_\nu$ is as in (\ref{eq:H_nu}), $\mathbf{S}_2$ is the matrix
defined through equation (\ref{eq:S_2_def}), and $\delta_n \to 0$.
Then, by a careful examination of the different terms in an
expansion of $H_\nu \mathbf{S}_2 \theta_\nu$, it is shown that an
upper bound on $\mathbb{E}\parallel H_\nu \mathbf{S}_2
\theta_\nu\parallel^2$ is asymptotically same as the RHS of
(\ref{eq:ASPCA_minimax_bound}). This is done in Section
\ref{subsec:aspca-H_nu_S_2_theta_nu}. Some results related to the
determination of correct asymptotic order of the terms in the
aforementioned expansion are given in Section
\ref{subsec:aspca-lemmas}. Before going into the detailed analysis,
it is necessary to fix some notation.

\subsection{Notation}

For any symmetric matrix $D$, $\lambda_k(D)$ will denote the $k$-th
largest eigenvalue of $D$. Frequently, the set $\{1,\ldots,N\}$ will
be divided into complementary sets $A$ and $B$. Here $A$ may refer
to the set of coordinates selected either in the first stage, or in
the second stage, or in a combination of both. $\mathbf{S}$ will be
partitioned as
\begin{equation}\label{eq:S_partition}
\mathbf{S} = \begin{bmatrix}
\mathbf{S}_{AA} & \mathbf{S}_{AB} \\
\mathbf{S}_{BA} & \mathbf{S}_{BB} \\
\end{bmatrix}
\end{equation}
where $\mathbf{S}_{AB}$ is the submatrix of $\mathbf{S}$ whose row
indices are from set $A$, and column indices are from set $B$. Any
$N \times 1$ vector $\mathbf{x}$ may similarly be partitioned as
$\mathbf{x} = (\mathbf{x}^T : \mathbf{y}^T)^T$. And for an $N \times
k$ matrix $\mathbf{Y}$, $\mathbf{Y}_A$ and $\mathbf{Y}_B$ will
denote the parts corresponding to rows with indices from set $A$ and
$B$, respectively. It should be clear, however, that no specific
order relation among these indices is assumed, and in fact the order
of the rows is unchanged in all of these situations. Expressions
like (\ref{eq:S_partition}) are just for convenience of writing.

\subsection{Bracketing relations}\label{sec:aspca-bracketing}

In this section the bracketing relationship is established. The
proof involves several parts. It essentially boils down to
probabilistic analysis of $1^o$ - $5^o$ of the ASPCA algorithm. This
is done in several stages. The coordinate selection step in $1^o$
and $2^o$ are jointly referred to as the \textit{first stage}, and
steps $3^o$, $4^o$ and $5^o$ are jointly referred to as the
\textit{second stage}.

\subsection{First stage coordinate selection}
\label{subsec:aspca-first_stage}

In this section  $1^o$, i.e., the first stage of the coordinate
selection scheme, is analyzed. Define
\begin{equation}\label{eq:zeta_k_def}
\zeta_k = \sum_{\nu=1}^M \lambda_\nu \theta_{\nu k}^2 ,
~~~k=1,\ldots,M.
\end{equation}
For $0 < a_- < 1 < a_+$, define
\begin{equation}\label{eq:SPCA_I_pm}
I_{1,n}^{\pm} = \{k : \zeta_k > a_{\mp} \gamma_1 \sqrt{\frac{\log(N
    \vee n)}{n}} \}.
\end{equation}
It is shown that $\widehat I_{1,n}$ satisfies the bracketing
relation (\ref{eq:bracketing_1}).

Let $\sigma_k^2 := \zeta_k + 1$. The selected coordinates are
\begin{equation}\label{eq:SPCA_select}
\widehat I_{1,n} = \{ k : \mathbf{S}_{kk} > 1 + \gamma_1
\sqrt{\frac{\log(N \vee n)}{n}} \} .
\end{equation}
Note that, $\mathbf{S}_{kk} \sim \sigma_k^2 \chi^2_{(n)}/n$. Then,
\begin{eqnarray}\label{eq:SPCA_FE}
\mathbb{P}(I_{1,n}^- \not\subset \widehat I_{1,n}) &=&
\mathbb{P}(\cup_{k \in I_{1,n}^-} \{ \mathbf{S}_{kk} \leq 1 +
\gamma_{1,n}\}) ~\leq~ \sum_{k \in I_{1,n}^-}
\mathbb{P}(\mathbf{S}_{kk} \leq 1 +
\gamma_{1,n})\nonumber\\
&\leq& \sum_{k \in I_{1,n}^-} \mathbb{P}(\frac{\mathbf{S}_{kk}}
{\sigma_k^2} \leq \frac{1+ \gamma_{1,n}}{1 + a_+
  \gamma_{1,n}})\nonumber\\
&\leq& |I_{1,n}^-| \mathbb{P}(\frac{\chi_{(n)}^2}{n} - 1\leq -
\frac{\gamma_{1,n}(a_+ - 1)}{1 + a_+  \gamma_{1,n}}), \qquad
(~\mbox{since},~
\mathbf{S}_{kk} \sim \sigma_k^2 \chi_{(n)}^2/n~) \nonumber\\
&\leq& |I_{1,n}^-| \exp\left(-\frac{n\gamma_{1,n}^2(a_+ - 1)^2}{4(1
+ a_+
    \gamma_{1,n})^2}\right), \qquad (\mbox{by}~
(\ref{eq:large_dev_chisq_2})~) \nonumber\\
&\leq& |I_{1,n}^-| (N \vee n)^{-(\gamma_1^2(a_+ - 1)^2/4)(1+o(1))}.
\end{eqnarray}
Similarly, if $n \geq 16$ then,
\begin{eqnarray}\label{eq:SPCA_FI}
\mathbb{P}(\widehat I_{1,n} \not\subset I_{1,n}^+) &=&
\mathbb{P}(\cup_{k
  \not\in I_{1,n}^+} \{ \mathbf{S}_{kk} > 1 + \gamma_{1,n}\}) ~\leq~
\sum_{k \not\in I_{1,n}^+} \mathbb{P}(\mathbf{S}_{kk} > 1 +
\gamma_{1,n})
\nonumber\\
&\leq& \sum_{k \not\in I_{1,n}^+}
\mathbb{P}(\frac{\mathbf{S}_{kk}}{\sigma_k^2} >
\frac{1+\gamma_{1,n}}{1+ a_- \gamma_{1,n}}) ~\leq~  N
\mathbb{P}(\frac{\chi_{(n)}^2}{n} - 1 > \frac{\gamma_{1,n}(1 -
a_-)}{1+
  a_- \gamma_{1,n}}) \nonumber\\
&\leq& N \frac{\sqrt{2}}{\gamma_1\sqrt{\log (N \vee n)}} \exp\left(-
\frac{n\gamma_{1,n}^2(1-a_-)^2}{4(1+a_-
    \gamma_{1,n})^2}\right), \qquad (\mbox{by}~
(\ref{eq:large_dev_chisq_3})~) \nonumber\\
&\leq& N (N \vee n)^{-(\gamma_1^2(1-a_-)^2/4)(1+o(1))}.
\end{eqnarray}
Combine (\ref{eq:SPCA_FE}) and (\ref{eq:SPCA_FI}) to get, as  $n \to
\infty$,
\begin{eqnarray}\label{eq:bracketing_1}
&&1-\mathbb{P}(I_{1,n}^- \subset \widehat I_{1,n} \subset I_{1,n}^+)
\nonumber\\
&\leq& |I_{1,n}^-| (N \vee n)^{-(\gamma_1^2(a_+ - 1)^2/4)(1+o(1))} +
 N (N \vee n)^{-(\gamma_1^2(1-a_-)^2/4)(1+o(1))} .
\end{eqnarray}
For future use, it is important to have an upper bound on the size
of the sets $I_{1,n}^{\pm}$. To this end, let $\mathbf{c} =
(c_1,\ldots,c_M)$ be such that $c_\nu > 0$ for all $\nu$ and
$\sum_{\nu=1}^M c_\nu^2 = 1$.
$$
I_{1,n}^{\pm} = \{ k \in \{1,\ldots,N\} : \sum_{\nu=1}^M \lambda_\nu
\theta_{\nu k}^2 > a_{\mp} \gamma_{1,n}\} \subset \bigcup_{\nu=1}^M
\{ k \in \{1,\ldots,N\} : |\theta_{\nu k}|
> c_\nu \sqrt{\frac{a_{\mp} \gamma_{1,n}}{\lambda_\nu}}\}.
$$
Since $\theta \in \Theta_q^M(C_1,\ldots,C_q)$, and $l^q(C)
\hookrightarrow wl^q(C)$, it follows from above that,
\begin{equation}\label{eq:J_1_pm}
|I_{1,n}^\pm| \leq J_{1,n}(\mathbf{c},\gamma_1,a_{\mp}) :=
a_\mp^{-q/2} \gamma_1^{-q/2} (\sum_{\nu=1}^M c_\nu^{-q}
\lambda_\nu^{q/2} C_\nu^q ) \frac{n^{q/4}}{(\log (N \vee n))^{q/4}}
.
\end{equation}
In fact, the upper bound is of the form
$J_{1,n}(\mathbf{c},\gamma_1,a_{\mp}) \wedge N$, since there are
altogether $N$ coordinates. Set $\mathbf{c} =
(M^{-1/2},\ldots,M^{-1/2})$, and denote the corresponding
$J_{1,n}(\mathbf{c},\gamma_1,a_\mp)$ by $J_{1,n}(\gamma_1,a_\mp)$.
Whenever there is no ambiguity about the choice of $\gamma_1$ and
$a_\mp$, $J_{1,n}(\gamma_1,a_\mp)$ will be denoted by $J_{1,n}^\pm$.
Notice that {\bf C1} and {\bf C2} imply that $J_{1,n}^+ \to \infty$
as $n \to \infty$. And {\bf C3} implies that
$\frac{J_{1,n}^+}{nh(\lambda_1)} \to 0$.

\vskip.1in\noindent{\bf Remark :}\label{rem:G_1_n_bound} From now
onwards, the set $\{I_{1,n}^- \subset \widehat I_{1,n} \subset
I_{1,n}^+\}$ will be denoted by $G_{1,n}$. Observe that $G_{1,n}$
depends on $\theta$. However, from (\ref{eq:bracketing_1}), it
follows that, if $\gamma_1 = 4$, $a_+ > 1+\frac{1}{\sqrt{2}}$ and $0
< a_- < 1-\frac{1}{\sqrt{2}}$, then there is an $\epsilon_0 > 0$ and
an $n_0\geq 1$, that depend on $a_+$ and $a_-$, such that for $n
\geq n_0$,
\begin{equation}\label{eq:G_1_n_bound}
\mathbb{P}(G_{1,n}^c) \leq (N \vee n)^{-1-\epsilon_0},
\end{equation}
uniformly in $\theta \in \Theta_q^M(C_1,\ldots,C_M)$.

\subsection{Eigen-analysis of $\mathbf{S}_{\widehat I_{1,n},\widehat
    I_{1,n}}$}\label{subsec:aspca-eigen_S_I_hat}

Throughout we follow the convention that $\langle \mathbf{e}_\nu,
\theta_{\nu,\widehat I_{1,n}} \rangle \geq 0$. Define
\begin{equation}\label{eq:S_1_def}
\widetilde{\mathbf{S}}_1 :=
\begin{bmatrix} \mathbf{S}_{\widehat I_{1,n} \widehat I_{1,n}} &   O \cr
                       O              &   O \cr
\end{bmatrix}
\qquad \mathbf{S}_1 :=
\begin{bmatrix} \mathbf{S}_{\widehat I_{1,n} \widehat I_{1,n}} &   O \cr
                       O              &   I \cr
\end{bmatrix} .
\end{equation}
Let $\widetilde{\mathbf{e}}_k$ be the eigenvector associated with
eigenvalue $\widehat \ell_k$ of $\widetilde{\mathbf{S}}_1$, for
$k=1,\ldots,m_1$, where $m_1 = (n \wedge |\widehat I_{1,n}|)$.
Eigenvalues of $\mathbf{S}_1$ belong to the set $\{\widehat
\ell_1,\ldots,\widehat \ell_{m_1}\} \cup \{1\}$; and the eigenvector
corresponding to the eigenvalue $\widehat \ell_k$ is
$\widetilde{\mathbf{e}}_k$, $1\leq k \leq m_1$. Note that, $\widehat
\ell_k$ is not necessarily the $k$-th largest eigenvalue of
$\mathbf{S}_1$. However, the analysis here will show that this
happens with very high probability for sufficiently large $n$.

Let $t_{1,n}^+ = 6(J_{1,n}^+/n \vee 1) \sqrt{\log(n \vee
J_{1,n}^+)/(n \vee J_{1,n}^+)}$. Define,
\begin{eqnarray}\label{eq:S_1_epsilons}
\varepsilon_{1,n} &=& \frac{2\sqrt{2}}{\lambda_1}\sum_{\nu=1}^M
\lambda_\nu
\sqrt{\frac{\log (n \vee J_{1,n}^+)}{n}}\nonumber\\
\varepsilon_{2,n} &=& \frac{2}{\lambda_1}\sum_{\nu=1}^M
\sqrt{\lambda_\nu}  \left(1+2\sqrt{2}\sqrt{\frac{\log (n \vee
J_{1,n}^+)} {n \wedge J_{1,n}^+}}\right) \sqrt{\frac{J_{1,n}^+}{n}}
\nonumber\\
\varepsilon_{3,n} &=& \frac{\sqrt{2}}{\lambda_1} \sum_{\nu\neq \nu'}
\sqrt{\lambda_\nu\lambda_{\nu'}}\sqrt{\frac{\log n}{n}}
\nonumber\\
\varepsilon_{4,n} &=&
\frac{1}{\lambda_1}\left(2\sqrt{\frac{J_{1,n}^+}{n}} +
\frac{J_{1,n}^+}{n} + \sqrt{2}t_{1,n}^+\right)
\nonumber\\
\varepsilon_{5,n} &=& c_q a_+^{1-q/2} \gamma_1^{1-q/2}
(\sum_{\nu=1}^M \lambda_\nu^{q/2} C_\nu^q) \frac{(\log(N \vee
n))^{1/2-q/4}} {\lambda_1 n^{1/2-q/4}}
\end{eqnarray}
where $c_q = \frac{2}{2-q}$. Observe that, under conditions {\bf
C1}-{\bf C3}, $\max_{1\leq j \leq 5} \varepsilon_{j,n} \to 0$ as $n
\to \infty$.

\noindent Set $A = \widehat I_{1,n}$, $A_+ = I_{1,n}^+$, $B =
\widehat I_{1,n}^c = \{1,\ldots,N\} \setminus \widehat I_{1,n}$, and
define
\begin{eqnarray} \label{eq:G_2_n}
G_{2,n} &=& \bigcap_{\nu=1}^M \{|\frac{\parallel v_\nu
\parallel^2}{n}-1| \leq 2\sqrt{2}\sqrt{\frac{\log n}{n}}\}
\nonumber\\
&& \bigcap~  \bigcap_{\nu=1}^M
 ~\{ \frac{\parallel  \mathbf{Z}_{A_+}
v_\nu \parallel}{n} \leq \left(1+ 2\sqrt{2}\sqrt{\frac{\log (n \vee
J_{1,n}^+)}{n \wedge J_{1,n}^+}}\right)
\sqrt{\frac{J_{1,n}^+}{n}} \}\nonumber\\
&& \bigcap~ \bigcap_{1\leq \nu < \nu'\leq M}
 ~\{|\frac{\langle v_\nu, v_{\nu'} \rangle}{n}| \leq
\sqrt{2} \sqrt{\frac{\log n}{n}}\},
\end{eqnarray}
and
\begin{equation}\label{eq:G_3_n}
G_{3,n} = \{ \frac{1}{n} \parallel \mathbf{Z}_{A_+}
\mathbf{Z}_{A_+}^T - I \parallel \leq 2\sqrt{\frac{J_{1,n}^+}{n}} +
\frac{J_{1,n}^+}{n} + \sqrt{2}t_{1,n}^+ \}.
\end{equation}
Then the following results hold.

\vskip.15in\noindent{\bf Lemma 3:}\label{lemma:S_1_eigen_deviation}
\textit{Under conditions {\bf C1}-{\bf C3},
\begin{eqnarray}
\bigcap_{j=1}^3 G_{j,n} &\subset& \{|\lambda_\nu(\mathbf{S}_1)  -
(1+\lambda_\nu)| \leq \lambda_1\sum_{j=1}^5 \varepsilon_{j,n}\},
\label{eq:S_1_eigen_deviation_set} \\
\mathbb{P}((G_{2,n} \cap G_{3,n})^c) &\leq& 3M (n \vee
J_{1,n}^+)^{-2} + M(M-1) n ^{-3+O(\frac{\log n}{n})} + 2 (n \vee
J_{1,n}^+)^{-2}. \label{eq:S_1_eigen_deviation}
\end{eqnarray}}

\vskip.15in\noindent{\bf Lemma 4:} \label{lemma:ell_hat_M_bound}
\textit{Let $\widetilde t_{1,n} = 6(|I_{1,n}^+|/n \vee
1)\sqrt{\log(n \vee |I_{1,n}^+|)/(n \vee |I_{1,n}^+|)}$. Under
conditions {\bf C1}-{\bf C3},
\begin{equation}\label{eq:ell_hat_M_bound}
\mathbb{P}(\widehat \ell_{M+1} > (1+\sqrt{\frac {|I_{1,n}^+|}{n}})^2
+ \sqrt{2} \widetilde t_{1,n}, ~\widehat{I}_{1,n} \subset I_{1,n}^+
) \leq 2(n \vee |I_{1,n}^+|)^{-2}.
\end{equation}}

\vskip.1in\noindent{\bf Remark :} \label{rem:eigen_S_AA_bound} Let
$G_{4,n} = \{\widehat \ell_{M+1} \leq (1+\sqrt{\frac
{|I_{1,n}^+|}{n}})^2 + \sqrt{2} \widetilde t_{1,n}\}$, where
$\widetilde t_{1,n}$ is as in Lemma 4.
Observe that $G_{4,n}$ depends on $\theta$; however,
$\mathbb{P}(G_{1,n} \cap G_{4,n}^c) \leq 2n^{-2}$ for all $\theta
\in \Theta_q^M(C_1,\ldots,C_M)$. It is easy to check that, under
{\bf C1}-{\bf C3},
\begin{equation}\label{eq:J_1_n_lambda_dominated}
2\sqrt{\frac{J_{1,n}^+}{n}} + \frac{J_{1,n}^+}{n} = o(\lambda_1)
~~~~\mbox{as}~~ n \to \infty.
\end{equation}
Therefore, from Lemma 3
and Lemma 4
it follows that, for sufficiently large $n$, uniformly in $\theta
\in \Theta_q^M(C_1,\ldots,C_M)$,
\begin{equation}\label{eq:ell_nu_hat_deviation}
\mathbb{P}(\max_{1\leq\nu \leq M} |\widehat \ell_\nu  -
(1+\lambda_\nu)| > \lambda_1\sum_{j=1}^5 \varepsilon_{j,n},
~G_{1,n}) \leq K_1(M)n^{-2},
\end{equation}
for some constant $K_1(M)$ that does not depend on $\theta$.

\subsection{Consistency of $\widehat M$}\label{subsec:aspca-consistency_M_hat}

\noindent{\bf Proposition 2:} \label{prop:M_hat_consistency}
\textit{Under conditions {\bf C1}-{\bf C3}, and with $\alpha_n$
defined through (\ref{eq:M_alpha_n}), $\widehat M$ is a consistent
estimator of $M$. In particular, if~ $\overline{\gamma}_1 = 9$,
$\gamma_1' = 3$, then there are constants $\overline{a}_+ > 1 >
\overline{a}_- > 0$, $1 > a' >0$, and an $n_{*0}$ such that for $n
\geq n_{*0}$, uniformly in $\theta \in \Theta_q^M(C_1,\ldots,C_M)$,
\begin{equation}\label{eq:M_hat_consistency}
\mathbb{P}(\widehat M  \neq M) \leq K_2(M) n^{-1-\epsilon_1},
\end{equation}
for some constants  $K_2(M) > 0$ and $\epsilon_1 :=
\epsilon_1(\overline{\gamma}_1, \gamma_1',\overline{a}_\pm, a')> 0$
independent of $\theta$.}

\subsection{Second stage coordinate selection}\label{subsec:aspca-second_stage}

Steps $4^0$ and $5^0$ of the ASPCA scheme are analyzed in this
subsection. For future reference, it is convenient to denote the
event $\bigcap_{j=1}^4 G_{j,n} \cap \{\widehat M = M\}$ by
$\overline{G}_{1,n}$. The ultimate goal of this section is to
establish (\ref{eq:second_stage_bracketing}). Throughout, it is
assumed that {\bf BA} and {\bf C1}-{\bf C3} are valid. Observe that,
by definition (see $4^o$ and $5^o$ of ASPCA scheme), $T_k =
\sum_{\mu=1}^M Q_{k\mu}^2$ if $k \not\in \widehat I_{1,n}$, and
define it to be zero otherwise.

\subsubsection{A preliminary bracketing relation}

First, define
\begin{equation}\label{eq:zeta_tilde_def}
\widetilde \zeta_k = \sum_{\nu=1}^M h(\lambda_\nu) \theta_{\nu k}^2,
~~~k=1,\ldots,N.
\end{equation}
Define, for $0 < \gamma_{2,-} < \gamma_2 < \gamma_{2,+}$,
\begin{equation}\label{eq:I_n_pm_def}
I_n^\pm = \{ k : \widetilde \zeta_k > \gamma_{2,\mp}^2 ~\frac{\log(N
  \vee n)}{n}\}.
\end{equation}
Observe that $\widetilde \zeta_k \geq \eta(\lambda_M) \zeta_k$. This
implies that, for some $n_{*1} \geq n_{*0} \vee n_{*0}'$, for all $n
\geq n_{*1}$, $I_{n,1}^+ \subset I_n^-$, uniformly in $\theta \in
\Theta_q^M(C_1,\ldots,C_M)$. Note that
\begin{equation*}
\mathbb{P}(\{I_n^- \subset \widehat I_{1,n} \cup \widehat I_{2,n}
\subset I_n^+\}^c, \overline{G}_{1,n}) \leq \mathbb{P}(I_n^-
\not\subset \widehat I_{1,n} \cup \widehat
I_{2,n},\overline{G}_{1,n}) + \mathbb{P}(\widehat I_{1,n} \cup
\widehat I_{2,n} \not \subset I_n^+,\overline{G}_{1,n}).
\end{equation*}
In the following, $D$ is a generic measurable set w.r.t. the
$\sigma$-algebra generated by $\mathbf{Z}$ and $v_1,\ldots,v_M$.
Then, for $n \geq n_{*1}$,
\begin{eqnarray}\label{eq:ASPCA_FI_inter}
\mathbb{P}(\widehat I_{1,n} \cup \widehat I_{2,n} \not \subset
I_n^+,~\overline{G}_{1,n} \cap D) &=& \mathbb{P}(\cup_{k \not\in
I_n^+} \{k \in \widehat I_{1,n} \cup
\widehat I_{2,n}\},~\overline{G}_{1,n}) \nonumber\\
= \mathbb{P}(\cup_{k \not\in I_n^+} \{k \in \widehat I_{2,n}\cap
\widehat I_{1,n}^c\}, ~\overline{G}_{1,n} \cap D) &\leq& \sum_{k
\not\in I_n^+} \mathbb{P}(k \in \widehat I_{2,n}\cap \widehat
I_{1,n}^c,
~\overline{G}_{1,n} \cap D) \nonumber\\
= \sum_{k \not\in I_n^+} \mathbb{P}(T_k > \gamma_{2,n}^2 ,
~\overline{G}_{1,n} \cap D), &&
\end{eqnarray}
where the last equality is from the inclusion $\widehat I_{1,n}
\subset I_{1,n}^+ \subset I_n^- \subset I_n^+$. Similarly,
\begin{eqnarray}\label{eq:ASPCA_FE_inter}
\mathbb{P}(I_n^- \not\subset \widehat I_{1,n} \cup \widehat I_{2,n},
~\overline{G}_{1,n} \cap D) &=& \mathbb{P}(\cup_{k \in I_n^-} \{k
\not\in \widehat
I_{1,n} \cup \widehat I_{2,n}\}, ~\overline{G}_{1,n} \cap D) \nonumber\\
= \mathbb{P}(\cup_{k \in I_n^- \setminus I_{1,n}^-} \{ k \in
\widehat I_{1,n}^c \cap \widehat I_{2,n}^c \}, ~\overline{G}_{1,n}
\cap D) &\leq& \sum_{k \in I_n^- \setminus I_{1,n}^-} \mathbb{P}(k
\not\in \widehat I_{1,n},k \not\in \widehat I_{2,n},
~\overline{G}_{1,n} \cap D)
\nonumber\\
= \sum_{k \in I_n^- \setminus I_{1,n}^-} \mathbb{P}(T_k \leq
\gamma_{2,n}^2, k \not\in \widehat I_{1,n}, ~\overline{G}_{1,n} \cap
D). &&
\end{eqnarray}

\subsubsection{Final bracketing relation}

It can be shown using some rather lengthy technical arguments
(provided in the technical note) that, given appropriate $\gamma_2$,
$\gamma_{2,+}$ and $\gamma_{2,-}$, for all sufficiently large $n$,
except on a set of negligible probability, uniformly in $\theta \in
\Theta_q^M(C_1,\ldots,C_M)$,
\begin{equation}\label{eq:T_k_comparison}
\begin{cases}
T_k < \gamma_{2,n}^2 & ~~if~~  k \not\in I_n^+, \\
T_k > \gamma_{2,n}^2 & ~~if~~  k \in I_n^- \setminus I_{1,n}^-. \\
\end{cases}
\end{equation}
Once (\ref{eq:T_k_comparison}) is established, it follows from
(\ref{eq:ASPCA_FI_inter}), (\ref{eq:ASPCA_FE_inter}),
and some probabilistic bounds (also given in the technical note)
that there exists $n_{*6}$  such that for all $n \geq n_{*6}$,
\begin{equation}\label{eq:second_stage_bracketing}
\mathbb{P}(I_n^- \subset \widehat I_{1,n} \cup\widehat I_{2,n}
\subset I_n^+,~\overline{G}_{1,n}) \geq 1 - K_6(M)
n^{-1-\epsilon_2(\kappa)},
\end{equation}
for some $K_6(M) > 0$ and $\epsilon_2(\kappa) > 0$. Moreover, the
bound (\ref{eq:second_stage_bracketing}) is uniform in $\theta \in
\Theta_q^M(C_1,\ldots,C_M)$.

\subsection{Second stage : perturbation analysis}\label{sec:aspca-perturbation}

The rest of this section deals with the part of the proof of Theorem
3
that involves analyzing the behavior of the submatrix of
$\mathbf{S}$ that corresponds to the set of selected coordinates. To
begin with, define $\widehat I_n := \widehat I_{1,n} \cup \widehat
I_{2,n}$, and $\overline{G}_{3,n} := \{I_n^- \subset \widehat I_n
\subset I_n^+\} \cap \overline{G}_{2,n}$.  Then define
\begin{equation}\label{eq:S_2_def}
\widetilde{\mathbf{S}}_2 =
\begin{bmatrix} \mathbf{S}_{\widehat I_n, \widehat I_n} &   O \cr
                       O              &   O \cr
\end{bmatrix}
\qquad \mathbf{S}_2 =
\begin{bmatrix} \mathbf{S}_{\widehat I_n, \widehat I_n} &   O \cr
                       O              &   I \cr
\end{bmatrix}  .
\end{equation}
In this section $A$ will denote the set $\widehat I_n$, $B =
\{1,\ldots,N\}\setminus A =: A^c$, $A_\pm = I_n^\pm$,
$\overline{A}_- = A_- \setminus A$, $B_- = \{1,\ldots,N\} \setminus
A_- =: A_-^c$. The first task before us is to derive an equivalent
of Lemma 3.
This is done in Section \ref{subsec:aspca-eigen_S_2}. The vector
$H_\nu \mathbf{S}_2 \theta_\nu$ is expanded, and then the important
terms are isolated in Section \ref{subsec:aspca-H_nu_S_2_theta_nu}.
Finally, the proof is completed in Section
\ref{subsec:aspca-proof_minimax}.

\subsection{Eigen-analysis of $\mathbf{S}_2$}\label{subsec:aspca-eigen_S_2}

$\widehat \theta_\nu$, $\nu=1,\ldots,M$ are the eigenvectors
corresponding to the  $M$ largest (in decreasing order) eigenvalues
of $\mathbf{S}_2$. As a convention $\langle \widehat \theta_\nu,
\theta_\nu \rangle \geq 0$ for all $\nu=1,\ldots,M$. Let the first
$M$ eigenvalues of $\widetilde{\mathbf{S}}_2$ be $\widetilde \ell_1
> \ldots > \widetilde \ell_M$. Then arguments similar to what are
used in Section \ref{subsec:aspca-eigen_S_I_hat} establishes the
following results.

On $\overline{G}_{3,n}$, for all $\mu=1,\ldots,M$,
\begin{equation}\label{eq:theta_mu_B_norm_bound_2}
\parallel \theta_{\mu,A_-^c}\parallel^2 \leq \overline\tau_{n,\mu}^2
:= c_q \gamma_{2,+}^{2-q} \frac{C_\mu^q (\log (n \vee
  N))^{1-q/2}}{(nh(\lambda_\mu))^{1-q/2}},
\end{equation}
and
\begin{equation}\label{eq:J_2_pm}
|I_n^\pm| \leq J_{2,n}^\pm := \gamma_{2,\mp}^{-q} M^{q/2}
(\sum_{\mu=1}^M h(\lambda_\mu)^{q/2} C_\mu^q) \left(\frac{n}{\log (n
\vee
    N)}\right)^{q/2}.
\end{equation}
Under {\bf C1}-{\bf C3}, as $n \to \infty$, for all
$\nu=1,\ldots,M$,
\begin{equation}\label{eq:J_2_pm_h_lambda_nu_bound}
\frac{J_{2,n}^\pm}{nh(\lambda_\nu)} \leq \gamma_{2,\mp}^{-q/2}
M^{q/2} \frac{\lambda_1^q}{\lambda_\nu^q} \frac{
(\sum_{\mu=1}^M(\frac {\lambda_\mu}{\lambda_1})^{q/2} C_\mu^q)(\log
(n \vee N))^{-q/2}} {(nh(\lambda_\nu))^{1-q/2}} \to 0;
\end{equation}
and $\overline\tau_n := \max_{1\leq \mu \leq M}\overline\tau_{n,\mu}
\to 0$. Again, check that $|I_n^\pm|$ is bounded by $N$, and
$\parallel \theta_{\mu,(I_n^-)^c}\parallel^2$ is bounded by
$\gamma_{2,+}^2 N\log (n\vee N) (nh(\lambda_\mu))^{-1}$. This
observation leads to the fact alluded to in Remark
\ref{rem:rate_comparison_OPCA}.

For $j=1,\ldots,4$, define $\overline\varepsilon_{j,n}$ as
$\varepsilon_{j,n}$ is defined in (\ref{eq:S_1_epsilons}), with
$J_{1,n}^+$ replaced by $J_{2,n}^+$. Then define
\begin{equation}\label{eq:S_2_epsilon_5}
\overline\varepsilon_{5,n} = c_q \gamma_{2,+}^{1-q/2}
\left[\sum_{\mu=1}^M
\left(\frac{\eta(\lambda_1)}{\eta(\lambda_\mu)}\right)^{1-q/2}
\left(\frac{\lambda_\mu}{\lambda_1}\right)^{q/2} C_\mu^q \right]
\left(\frac{\log (n \vee N)}{nh(\lambda_1)}\right)^{1-q/2}.
\end{equation}
It follows that $\max_{1\leq j \leq 5} \overline\varepsilon_{j,n}
\to 0$ as $n \to \infty$. Define
\begin{equation}\label{eq:Delta_n_nu_bar}
\overline \Delta_{n,\nu} = \frac{\lambda_1}{\max\{\lambda_{\nu-1} -
  \lambda_\nu, \lambda_\nu - \lambda_{\nu+1}\}}
\left[\sum_{j=1}^5 \overline\varepsilon_{j,n} +
  \sqrt{\sum_{\mu=1}^M
    \frac{\lambda_\mu}{\lambda_1}}\sqrt{\overline\varepsilon_{5,n}} \right],
\end{equation}
and $\overline \Delta_n = \max_{1\leq \nu \leq M} \overline
\Delta_{n,\nu}$. A result that summarizes the behavior of the first
$M$ eigenvalues of $\mathbf{S}_2$ can now be stated.

\vskip.15in\noindent{\bf Proposition
3:}\label{prop:S_2_eigen_bounds} \textit{There is a measurable set
$\overline G_{4,n} \subset \overline G_{3,n}$, and an integer
$n_{*7} \geq n_{*6}$, such that, for all $n \geq n_{*7}$ the
following relations hold, uniformly in $\theta \in
\Theta_q^M(C_1,\ldots,C_M)$.
\begin{eqnarray}
\overline{G}_{4,n} &\subset& \bigcap_{\nu=1}^M \{\widetilde \ell_\nu
= \lambda_\nu(\mathbf{S}_2)~\mbox{and}~|\widetilde \ell_\nu -
(1+\lambda_\nu)| \leq \lambda_1
\sum_{j=1}^5\overline\varepsilon_{j,n}\},\label{eq:S_2_eigen_deviation_set}\\
\overline{G}_{4,n} &\subset& \{\parallel \mathbf{S}_2 - \Sigma
\parallel \leq \lambda_1 (\sum_{j=1}^5 \overline\varepsilon_{j,n} +
  \sqrt{\sum_{\mu=1}^M
    \frac{\lambda_\mu}{\lambda_1}}\sqrt{\overline\varepsilon_{5,n}})\}
\label{eq:S_2_Sigma_diff_bound}\\
1-\mathbb{P}(\overline G_{4,n}) &\leq& K_7(M)n^{-1-\epsilon_3},
 \label{eq:S_2_eigen_bounds}
\end{eqnarray}
for some constants $K_7(M) > 0$ and $\epsilon_3 > 0$. $\epsilon_3$
depends of $\gamma_1$, $\overline\gamma_1$, $\gamma_1'$, $a_\pm$,
$\gamma_2$, $\gamma_{2,\pm}$, and $\kappa$.}

\vskip.15in At this point it is useful to define a quantity that
will play an important role in the analysis in Section
\ref{subsec:aspca-H_nu_S_2_theta_nu}. Define,
\begin{equation}\label{eq:vartheta_mu_def}
\vartheta_{n,\mu}^2 = \overline{\tau}_{n,\mu}^2 +
\frac{J_{2,n}^+}{nh(\lambda_\mu)} + \sum_{\mu'\neq\mu}
\frac{1}{ng(\lambda_{\mu'},\lambda_\mu)}, ~~~~\mu=1,\ldots,M.
\end{equation}
Then define $\vartheta_n = \max_{1\leq \mu \leq M}
\vartheta_{n,\mu}$ and observe that, under {\bf C1}-{\bf C2},
$\vartheta_n \to 0$ as $n \to \infty$.

We argue that,  for $n \geq n_{*8}$, say, on a set
$\overline{G}_{5,n}$ with probability approaching 1 sufficiently
fast,
\begin{equation}\label{eq:theta_hat_nu_theta_nu_loss_diff}
L(\widehat \theta_\nu,\theta_\nu) \leq \parallel H_\nu \mathbf{S}_2
\theta_\nu \parallel^2 (1+ \overline\delta_{n,N,\nu}),
\end{equation}
where $\overline\delta_{n,N,\nu} = o(1)$. Therefore, it remains to
show that, $\mathbb{E}\parallel H_\nu \mathbf{S}_2 \theta_\nu
\parallel^2 \mathbf{1}_{\overline{G}_{5,n}}$ is bounded by the
quantity appearing on the RHS of (\ref{eq:ASPCA_minimax_bound}).

\subsection{Analysis of $H_\nu \mathbf{S}_2
  \theta_\nu$}\label{subsec:aspca-H_nu_S_2_theta_nu}

In this section, as in Section \ref{subsec:aspca-proof_minimax},
$\nu$ is going to be a fixed index in $\{1,\ldots,M\}$.
Before an analysis of $H_\nu \mathbf{S}_2 \theta_\nu$ is carried
out, a few important facts are stated below. Here $C$ is any subset
of $\{1,\ldots,N\}$ satisfying $A_- \subset C$.
\begin{eqnarray}
|\delta_{\mu\nu} - \langle \theta_{\mu,C} ,\theta_{\nu,C}\rangle|
&=& |\langle \theta_{\mu,C^c} \theta_{\nu,C^c}\rangle| \leq
\overline\tau_{n,\mu} \overline\tau_{n,\nu} \leq (\vartheta_{n,\mu}
\vee \vartheta_{n,\mu})\vartheta_n,
\label{eq:theta_mu_theta_nu_inner_bound}\\
\max_{1\leq \mu\leq M}\frac{|I_n^\pm|}{n h(\lambda_\mu)} &\leq&
\frac{h(\lambda_\nu)}{h(\lambda_M)} \vartheta_{n,\nu}^2.
\label{eq:I_n_pm_varepsilon_bound}
\end{eqnarray}
Further,
\begin{equation}\label{eq:tau_n_bar_log_n_vartheta}
\max_{1\leq \mu,\mu' \leq M} \parallel \theta_{\mu,C}\parallel
\sqrt{\frac{\log n}{nh(\lambda_{\mu'})}} \leq \overline\tau_n
\sqrt{\frac{\log n}{nh(\lambda_M)}} = O(\frac{\log
n\sqrt{J_{2,n}^+}}{nh(\lambda_\nu)})=o(\vartheta_{n,\nu}),
\end{equation}
which follows from {\bf C1}, {\bf C2},
(\ref{eq:theta_mu_B_norm_bound_2}),
(\ref{eq:J_2_pm_h_lambda_nu_bound}), and (\ref{eq:vartheta_mu_def}).

Next, observe that $H_\nu\theta_\nu = 0$ implies that
\begin{equation}\label{eq:H_S_theta_nu}
H_\nu \mathbf{S}_2 \theta_\nu = H_\nu(\mathbf{S}_2 - \Sigma)
\theta_\nu = \begin{bmatrix} H_{\nu,AA} (\mathbf{S}_{AA}
-I)\theta_{\nu,A}  \cr H_{\nu,BA} (\mathbf{S}_{AA} -I)\theta_{\nu,A}
\cr
\end{bmatrix}
= \Psi, ~~\mbox{say}.
\end{equation}
Then $\Psi_A$ and $\Psi_B$ have the general form, for $C= A, B$,
\begin{eqnarray}\label{eq:Psi_C}
\Psi_C &=& \sum_{\mu=1}^M \frac{\parallel
  v_\mu \parallel^2}{n} \lambda_\mu
\langle \theta_{\mu,A},\theta_{\nu,A}\rangle H_{\nu,CA}
\theta_{\mu,A} + \sum_{\mu=1}^M \sqrt{\lambda_\mu} \frac{1}{n}
\langle \mathbf{Z}_A v_\mu, \theta_{\nu,A}
\rangle H_{\nu,CA} \theta_{\mu,A} \nonumber\\
&& + \sum_{\mu=1}^M \sqrt{\lambda_\mu} \langle
\theta_{\mu,A},\theta_{\nu,A} \rangle H_{\nu,CA} \frac{1}{n}
\mathbf{Z}_A v_\mu + \sum_{\mu \neq \mu'} \frac{\langle v_\mu,
  v_{\mu'}\rangle}{n} \sqrt{\lambda_\mu \lambda_{\mu'}}
\langle \theta_{\mu',A}, \theta_{\nu,A}\rangle H_{\nu,CA}
\theta_{\mu,A}
\nonumber\\
&& + H_{\nu,CA}(\frac{1}{n} \mathbf{Z}_A
\mathbf{Z}_A^T-I)\theta_{\nu,A}.
\end{eqnarray}
When $C$ is either $A$ or $B$, and $\delta_{CA}$ is 1 or 0 according
as whether $C=A$ or not,
\begin{eqnarray}
H_{\nu,CA} \theta_{\mu,A} &=& \sum_{\nu'\neq
  \nu}\frac{1}{\lambda_{\nu'}-\lambda_\nu} \langle
\theta_{\nu',A},\theta_{\mu,A}\rangle \theta_{\nu',C} +
\frac{1}{\lambda_\nu} \sum_{\nu'\neq \mu}^M \langle
\theta_{\nu',A},\theta_{\mu,A}\rangle \theta_{\nu',C} \nonumber\\
&&~~~~~~~~~~~~~~~~~~ -\frac{1}{\lambda_\nu}(\delta_{CA}-
\parallel \theta_{\mu,A}\parallel^2)
\theta_{\mu,C}; \label{eq:H_nu_CA_theta_mu}\\
H_{\nu,CA} \mathbf{Z}_Av_\mu &=& \sum_{\nu'\neq \nu}
\frac{1}{\lambda_{\nu'}-\lambda_\nu} \langle \mathbf{Z}_Av_\mu,
\theta_{\nu',A}\rangle \theta_{\nu',C} \nonumber\\
&& ~~~~~~~~~~~~~~~~-\frac{1}{\lambda_\nu}(\delta_{CA} I -
\sum_{\nu'=1}^M\theta_{\nu',C}\theta_{\nu',A})\mathbf{Z}_Av_\mu;
\label{eq:H_nu_CA_Z_A_v_mu}
\end{eqnarray}
\begin{eqnarray}\label{eq:H_nu_CA_Z_A_Z_A_theta_mu}
H_{\nu,CA}(\frac{1}{n} \mathbf{Z}_A \mathbf{Z}_A^T-I)\theta_{\nu,A}
&=& \sum_{\nu'\neq \nu} \frac{1}{\lambda_{\nu'}-\lambda_\nu}
\left(\frac{1}{n}\langle \mathbf{Z}_A^T \theta_{A,\nu'},
\mathbf{Z}_A^T \theta_{A,\nu}\rangle -\langle \theta_{\nu',A},
\theta_{\nu,A}\rangle \right) \theta_{\nu',C} \nonumber\\
&&  -\frac{1}{\lambda_\nu}(\delta_{CA} I -
\sum_{\nu'=1}^M\theta_{\nu',C}\theta_{\nu',A})(\frac{1}{n}
\mathbf{Z}_A \mathbf{Z}_A^T-I)\theta_{\nu,A}.
\end{eqnarray}
A further expansion of terms $\Psi_A$ and $\Psi_B$ can be computed,
but at this point it is beneficial to isolate the important terms in
the expansion. Accordingly, use Lemmas 9 - 13,
together with
(\ref{eq:theta_mu_theta_nu_inner_bound}),
(\ref{eq:I_n_pm_varepsilon_bound}) and
(\ref{eq:tau_n_bar_log_n_vartheta}) to deduce that, there is a
measurable set $\overline{G}_{5,n} \subset \overline{G}_{4,n}$,
constants $K_8(M) > 0$, $\epsilon_4 > 0$ and an $n_{*8} \geq n_{*7}$
such that, $1-\mathbb{P}(\overline{G}_{5,n}) \leq
K_8(M)n^{-1-\epsilon_4}$, for $n \geq n_{*8}$, and
\begin{equation}\label{eq:Psi_expand}
\Psi = \Psi_0 + \Psi_I + \Psi_{II} + \Psi_{III}  + \Psi_{IV} +
\Psi_{rem},
\end{equation}
where $\parallel \Psi_{rem}\parallel \leq b_n \vartheta_{n,\nu}$,
with $b_n = o(1)$, and the other elements are described below.

\begin{equation}\label{eq:Psi_0}
\Psi_{0,A} = 0 ~~~\mbox{and}~~~ \Psi_{0,B} = \theta_{\nu,B}.
\end{equation}
$\Psi_I = \sum_{\mu\neq \nu}^M w_{\mu\nu} \theta_\mu$ where
$w_{\mu\nu}$ equals
\begin{eqnarray}\label{eq:Psi_I}
&& \hskip-.1in \frac{\sqrt{\lambda_\mu}}{\lambda_\mu - \lambda_\nu}
\frac{1}{n} \langle \mathbf{Z}_{A_-} v_\mu, \theta_{\nu,A_-} \rangle
+ \frac{\sqrt{\lambda_\nu}}{\lambda_\mu - \lambda_\nu} \frac{1}{n}
\langle \mathbf{Z}_{A_-} v_\nu, \theta_{\mu,A_-} \rangle
\nonumber\\
&& \hskip-.1in + \frac{\sqrt{\lambda_\mu\lambda_\nu}}{\lambda_\mu -
\lambda_\nu} \frac{\langle v_\mu, v_\nu\rangle}{n} +
\frac{1}{\lambda_\mu -
  \lambda_\nu} \left(\frac{1}{n}\langle
\mathbf{Z}_{A_-}^T \theta_{A_-,\mu}, \mathbf{Z}_{A_-}^T
\theta_{A_-,\nu}\rangle -\langle \theta_{\mu,A_-},
\theta_{\nu,A_-}\rangle \right)
\end{eqnarray}
\begin{equation}\label{eq:Psi_II}
\Psi_{II} = -\frac{1}{\lambda_\nu} (I - \sum_{\mu=1}^M \theta_\mu
\theta_\mu^T) (\frac{1}{n}
\widetilde{\mathbf{Z}}\widetilde{\mathbf{Z}}^T - \Xi) \theta_\nu,
\end{equation}
where $\widetilde{\mathbf{Z}}_{A_-} = \mathbf{Z}_{A_-}$ and
$\widetilde{\mathbf{Z}}_{A_-^c} = O$, i.e. a matrix whose entries
are all 0; and $\Xi$ is a $N\times N$ matrix whose $(A_-,A_-)$ block
is identity and the rest are all zero.
\begin{equation}\label{eq:Psi_III}
\Psi_{III} = -\frac{1}{\sqrt{\lambda_\nu}} (I - \sum_{\mu=1}^M
\theta_\mu \theta_\mu^T) \frac{1}{n} \widetilde{\mathbf{Z}}v_\nu.
\end{equation}
$\Psi_{IV}$ is such that $\Psi_{IV,A_-} = 0$, $\Psi_{IV,B} = 0$, and
\begin{equation}\label{eq:Psi_IV}
\Psi_{IV,\overline{A}_-} = -\frac{1}{n}\mathbf{Z}_{\overline{A}_-}
\left(\frac{1}{\sqrt{\lambda_\nu}} v_\nu + \frac{1}{\lambda_\nu}
\mathbf{Z}_{A_-}^T \theta_{\nu,A_-}\right).
\end{equation}

\subsection{Completion of the proof of Theorem 3}
\label{subsec:aspca-proof_minimax}

Suppose without loss of generality that $n_{*8}$ in Section
\ref{subsec:aspca-H_nu_S_2_theta_nu} is large enough so that
$\overline\Delta_n < \frac{\sqrt{5}-1}{4}$. Since on
$\overline{G}_{5,n}$, $\parallel \mathbf{S}_2 - \Sigma \parallel
\leq \min\{\lambda_\nu - \lambda_{\nu+1},
\lambda_{\nu-1}-\lambda_\nu\} \overline{\Delta}_n$, where $\lambda_0
= \infty$ and $\lambda_{M+1} = 0$, argue that, by Lemma 6,
for $n \geq n_{*8}$, on
$\overline{G}_{5,n}$,
\begin{equation}\label{eq:theta_hat_nu_theta_nu_loss_diff}
L(\widehat \theta_\nu,\theta_\nu) \leq \parallel H_\nu \mathbf{S}_2
\theta_\nu \parallel^2 (1+ \overline\delta_{n,N,\nu}),
\end{equation}
where $\overline\delta_{n,N,\nu} = o(1)$. Therefore, it remains to
show that, $\mathbb{E}\parallel H_\nu \mathbf{S}_2 \theta_\nu
\parallel^2 \mathbf{1}_{\overline{G}_{5,n}}$ is bounded by the
quantity appearing on the RHS of (\ref{eq:ASPCA_minimax_bound}). In
view of the fact that, this upper bound is within a constant
multiple of $\vartheta_{n,\nu}^2$, and $\parallel
\Psi_{rem}\parallel = o(\vartheta_{n,\nu})$ on $\overline{G}_{5,n}$,
it is enough that the same bound holds for $\mathbb{E}\parallel \Psi
- \Psi_{rem}\parallel^2 \mathbf{1}_{\overline{G}_{5,n}}$.

Observe that $\Psi_{I}$, $\Psi_{II}$, and $\Psi_{III}$ are mutually
uncorrelated vectors. Also, by (\ref{eq:theta_mu_B_norm_bound_2}),
$\mathbb{E}\parallel
\Psi_0\parallel^2\mathbf{1}_{\overline{G}_{5,n}} \leq
\overline\tau_{n,\nu}^2$. Therefore,
\begin{eqnarray}\label{eq:Psi_expect_1}
&&\mathbb{E}\parallel \Psi - \Psi_{rem}\parallel^2
\mathbf{1}_{\overline{G}_{5,n}}
\nonumber\\
&\leq& \hskip-.1in \overline\tau_{n,\nu}^2 + \mathbb{E}\parallel
\Psi_I\parallel^2 + \mathbb{E}\parallel \Psi_{II}\parallel^2 +
\mathbb{E} \parallel \Psi_{III}\parallel^2 + \mathbb{E}\parallel
\Psi_{IV}\parallel^2\mathbf{1}_{\overline{G}_{5,n}} \nonumber\\
&& \hskip-.1in + 2\mathbb{E}|\langle \Psi_0, \Psi_{I} + \Psi_{II} +
\Psi_{III}\rangle| \mathbf{1}_{\overline{G}_{5,n}} +
2\mathbb{E}|\langle \Psi_{IV}, \Psi_0 + \Psi_{I} + \Psi_{II} +
\Psi_{III}\rangle| \mathbf{1}_{\overline{G}_{5,n}}
\end{eqnarray}
Observe that, $\Psi_{II,A_-} = -\frac{1}{\lambda_\nu} (I -
\sum_{\mu=1}^M \theta_{\mu,A_-} \theta_{\mu,A_-}^T) (\frac{1}{n}
\mathbf{Z}_{A_-}\mathbf{Z}_{A_-}-I)\theta_{\nu,A_-}$,
$$
\Psi_{II,A_-^c} = \frac{1}{\lambda_\nu}\sum_{\mu=1}^M
\left(\frac{1}{n} \langle \mathbf{Z}_{A_-}^T
\theta_{\mu,A_-},\mathbf{Z}_{A_-}^T \theta_{\nu,A_-}\rangle -\langle
\theta_{\mu,A_-},\theta_{\nu,A_-}\rangle \right) \theta_{\mu,A_-^c},
$$
and
$$
\Psi_{III,A_-} = -\frac{1}{\sqrt{\lambda_\nu}} (I - \sum_{\mu=1}^M
\theta_{\mu,A_-} \theta_{\mu,A_-}^T) \frac{1}{n}
\mathbf{Z}_{A_-}v_\nu, ~\Psi_{III,A_-^c} =
\frac{1}{\sqrt{\lambda_\nu}} \sum_{\mu=1}^M \frac{1}{n} \langle
\mathbf{Z}_{A_-}v_\nu, \theta_{\mu,A_-}\rangle \theta_{\mu,A_-^c}.
$$
Thus, by a further application of Lemmas 9-12,
it can be checked that, there is an integer $n_{*9}\geq n_{*8}$, and
an event $\overline{G}_{6,n} \subset \overline{G}_{5,n}$ such that,
for $n \geq n_{*9}$,  on $\overline{G}_{6,n}$,
\begin{equation}\label{eq:Psi_cross_prod_bound}
|\langle\Psi_0,\Psi_{I} + \Psi_{II} + \Psi_{III}\rangle|+ |\langle
\Psi_{IV},\Psi_0 + \Psi_{I} + \Psi_{II} + \Psi_{III}\rangle| \leq
b_n'\vartheta_{n,\nu}^2,
\end{equation}
with $b_n' = o(1)$; and $\mathbb{P}(\overline{G}_{6,n}^c \cap
\overline{G}_{5,n}) \leq n^{-2(1+\epsilon_5)}$, for some constants
$K_9(M) > 0$ and $\epsilon_5
> 0$.
On $\overline{G}_{5,n}$,
$$
\parallel \Psi_{IV}\parallel^2 \leq \frac{1}{n^2} \parallel
\mathbf{Z}_{A_{+/-}} \left(\frac{1}{\sqrt{\lambda_\nu}} v_\nu +
\frac{1}{\lambda_\nu} \mathbf{Z}_{A_-}^T
\theta_{\nu,A_-}\right)\parallel^2,
$$
where $A_{+/-} := A_+\setminus A_-$, and the (unrestricted)
expectation of the random variable appearing in the upper bound is
bounded by $\frac{|I_n^+| - |I_n^-|}{nh(\lambda_\nu)}$. From this,
and some expectation computations similar to those in Section
\ref{subsec:highdpca-OPCA_risk_bound}, deduce that,
\begin{equation}\label{eq:Psi_expect_2}
\overline\tau_{n,\nu}^2 + \mathbb{E}\parallel \Psi_I\parallel^2 +
\mathbb{E}\parallel \Psi_{II}\parallel^2 + \mathbb{E} \parallel
\Psi_{III}\parallel^2 + \mathbb{E} \parallel
\Psi_{IV}\parallel^2\mathbf{1}_{\overline{G}_{5,n}} \leq
\vartheta_{n,\nu}^2 (1+o(1)).
\end{equation}
Finally, express the event $\overline{G}_{5,n}$ as (disjoint) union
of $\overline{G}_{5,n}\cap \overline{G}_{6,n}$ and
$\overline{G}_{5,n}\cap \overline{G}_{6,n}^c$; apply the bound
(\ref{eq:Psi_cross_prod_bound}) for the first set, and use
Cauchy-Schwartz inequality for the second set, to conclude that,
\begin{eqnarray}\label{eq:Psi_expect_3}
&& \mathbb{E}|\langle \Psi_0, \Psi_{I} + \Psi_{II} +
\Psi_{III}\rangle| \mathbf{1}_{\overline{G}_{5,n}} +
\mathbb{E}|\langle \Psi_{IV}, \Psi_0 + \Psi_{I} + \Psi_{II} +
\Psi_{III}\rangle|
\mathbf{1}_{\overline{G}_{5,n}} \nonumber\\
&\leq& b_n'\vartheta_{n,\nu}^2 + 2\sqrt{K_9(M)}n^{-1-\epsilon_5}
\vartheta_n = o(\vartheta_{n,\nu}^2).
\end{eqnarray}
Combine (\ref{eq:Psi_expect_1}), (\ref{eq:Psi_expect_2}) and
(\ref{eq:Psi_expect_3}) to complete the proof.

\section{Appendix}

Some results that are needed to prove the three theorems are
presented here.

\subsection{Deviation of extreme eigenvalues}
\label{subsec:highdpca-eigen_deviation_bound}

The goal is to provide a probabilistic bound for deviations of
$\parallel \frac{1}{n} \mathbf{Z}\mathbf{Z}^T - I\parallel$. This is
achieved through the following lemma.

\vskip.1in\noindent{\bf Lemma 5:}
\label{lemma:eigen_deviation_bound} \textit{Let $t_n = 6
(\frac{N}{n} \vee 1) \sqrt{\frac{\log (n \vee N)}{n \vee N}}$. Then,
for any $c
>0$, there exists $n_c \geq 1$ such that, for all $n \geq n_c$,
\begin{equation}\label{eq:eigen_deviation_bound}
\mathbb{P}\left(\parallel \frac{1}{n} \mathbf{Z}\mathbf{Z}^T -
I\parallel > 2\sqrt{\frac{N}{n}} + \frac{N}{n} + c t_n \right) \leq
2(n \vee N)^{-c^2} .
\end{equation}}

\noindent {\bf Proof :} By definition,
\begin{equation*}
\parallel \frac{1}{n} \mathbf{Z}\mathbf{Z}^T - I\parallel
= \max \{ \lambda_1(\frac{1}{n} \mathbf{Z}\mathbf{Z}^T) - 1, 1 -
\lambda_N (\frac{1}{n} \mathbf{Z}\mathbf{Z}^T) \}.
\end{equation*}
From Proposition 4
(due to Davidson and Szarek (2001)), and its consequence, Corollary
2,
given below, it follows that,
\begin{eqnarray}\label{eq:eigen_deviation_bound_1}
&& \mathbb{P}\left(\parallel \frac{1}{n} \mathbf{Z}\mathbf{Z}^T
- I\parallel > 2\sqrt{\frac{N}{n}} + \frac{N}{n} + c t_n \right) \nonumber\\
&\leq& \exp\left(-\frac{nc^2 t_n^2}{8(c t_n +
(1+\sqrt{N/n})^2)}\right) + \exp\left(-\frac{nc^2 t_n^2}{8(c t_n +
(1-\sqrt{N/n})^2)}\right).
\end{eqnarray}
First suppose that $n \geq N$. Then for $n$ large enough, $ct_n <
\frac{1}{2}$, so that
$$
\frac{nc^2 t_n^2}{8(c t_n + (1+\sqrt{N/n})^2)} \geq \frac{n c^2
  t_n^2}{36} , \quad \mbox{and} \quad
\frac{nc^2 t_n^2}{8(c t_n + (1-\sqrt{N/n})^2)} \geq \frac{nc^2
  t_n^2}{12}.
$$
Since in this case $nt_n^2 = 36 \log n$,
(\ref{eq:eigen_deviation_bound}) follows from
(\ref{eq:eigen_deviation_bound_1}). If $N > n$, then
$\lambda_N(\frac{1}{n} \mathbf{Z}\mathbf{Z}^T) = 0$, and
$$
\frac{nc^2 t_n^2}{8(c t_n + (1\pm \sqrt{N/n})^2)} = \frac{N c^2
(\frac{n}{N} t_n)^2 } {8(c \frac{n}{N} t_n + (1\pm \sqrt{n/N})^2)}
~,
$$
and therefore, (\ref{eq:eigen_deviation_bound}) follows if the roles
of $n$ and $N$ are reversed.

\vskip.15in\noindent{\bf Proposition 4:}
\label{prop:extreme_singval_concen} \textit{Let $Z$ be a $p \times
q$ matrix of i.i.d. $N(0,1)$ entries with $p \leq q$. Let
$s_{max}(Z)$ and $s_{min}(Z)$  denote the largest and the smallest
singular value of $Z$, respectively. Then,
\begin{eqnarray}
\mathbb{P}(s_{max}(\frac{1}{\sqrt{q}} Z) > 1+ \sqrt{p/q} + t) &\leq&
e^{-qt^2/2} ,
\label{eq:extreme_singular_1}\\
\mathbb{P}(s_{min}(\frac{1}{\sqrt{q}}Z) < 1- \sqrt{p/q} - t) &\leq&
e^{-qt^2/2} . \label{eq:extreme_singular_2}
\end{eqnarray}}

\vskip.15in\noindent{\bf Corollary 2:}
\label{cor:extreme_eval_concen} \textit{Let $\mathbf{S} =
\frac{1}{q} Z Z^T$ where $Z$ is as in Proposition 4,
with $p \leq q$. Let $m_1(p,q) :
= (1+\sqrt{\frac{p}{q}})^2$ and $m_p(p,q) :=
(1-\sqrt{\frac{p}{q}})^2$. Let $\lambda_1(\mathbf{S})$ and
$\lambda_p(\mathbf{S})$ denote the largest and the smallest
eigenvalues of $\mathbf{S}$. Then, for $t > 0$,
\begin{eqnarray}\label{eq:large_eigen_concentration}
\mathbb{P}(\lambda_1(\mathbf{S}) - m_1(p,q) > t) &\leq&
\exp\left(-\frac{q}{2}(\sqrt{t+m_1(p,q)} - \sqrt{m_1(p,q)})^2\right)
\nonumber\\
&\leq& \exp\left(-\frac{qt^2}{8(t+m_1(p,q))}\right),
\end{eqnarray}
and
\begin{eqnarray}\label{eq:small_eigen_concentration}
\mathbb{P}(\lambda_p(\mathbf{S}) - m_p(p,q) < -t) &\leq&
\exp\left(-\frac{q}{2}(\sqrt{t+m_p(p,q)} - \sqrt{m_p(p,q)})^2\right)
\nonumber\\
&\leq& \exp\left(-\frac{qt^2}{8(t+m_p(p,q))}\right).
\end{eqnarray}}

\subsection{Perturbation of
  eigen-structure}\label{subsec:highdpca-perturb_eigen}

The following lemma is most convenient for the risk analysis of
estimators of $\theta_\nu$. Several variants of this lemma appear in
the literature (Kneip and Utikal (2001), Tyler (1983), Tony Cai and
Hall (2005)) and most of them implicitly use the approach proposed
by Kato (1980).

\vskip.1in\noindent {\bf Lemma 6:} \label{lemma:evec_perturb_bound}
\textit{For some $T \in \mathbb{N}$, let $A$ and $B$ be two
symmetric $T \times T$ matrices. Let the eigenvalues of matrix $A$
be denoted by $\lambda_1(A) \geq \ldots \geq \lambda_T(A)$. Set
$\lambda_0(A) = \infty$ and $\lambda_{T+1}(A) = - \infty$. For any
$r \in\{1,\ldots,T\}$, if $\lambda_r(A)$ is a unique eigenvalue of
$A$, i.e., if $\lambda_{r-1}(A) > \lambda_r(A) > \lambda_{r+1}(A)$,
then denoting by $\mathbf{p}_r$ the eigenvector associated with the
$r$-th eigenvalue,
\begin{equation}\label{eq:eig_perturb_first}
\mathbf{p}_r(A+B) - \mbox{sign}(\mathbf{p}_r(A+B)^T\mathbf{p}_r(A))
\mathbf{p}_r(A) = - H_r(A) B \mathbf{p}_r(A) + R_r
\end{equation}
where $H_r(A) := \sum_{s \neq r} \frac{1}{\lambda_s(A) -
\lambda_r(A)} P_{{\cal E}_s}(A)$ and $P_{{\cal E}_s}(A)$ denotes the
projection matrix onto the eigenspace ${\cal E}_s$ corresponding to
eigenvalue $\lambda_s(A)$ (possibly multi-dimensional). Define
$\Delta_r$ and  $\overline{\Delta}_r$ as
\begin{eqnarray}
\Delta_r &:=& \frac{1}{2} [\parallel H_r(A) B \parallel +
|\lambda_r(A+B)
- \lambda_r(A)| \parallel H_r(A) \parallel] \label{eq:eigen_Delta_r} \\
\overline{\Delta}_r &=& \frac{\parallel B\parallel} {\min_{1\leq j
\neq r \leq T} |\lambda_j(A) -
\lambda_r(A)|}~.\label{eq:eigen_Delta_bar_r}
\end{eqnarray}
Then, the residual term $R$ can be bounded by
\begin{equation}\label{eq:eigenvec_error}
\parallel R_r \parallel \leq \min\{ 10\overline{\Delta}_r^2, ~\parallel H_r(A)
B \mathbf{p}_r(A)
\parallel \left[\frac{2\Delta_r(1+2\Delta_r)}{1-2\Delta_r(1+2\Delta_r)} +
\frac{\parallel H_r(A) B \mathbf{p}_r(A)\parallel}
{(1-2\Delta_r(1+2\Delta_r))^2}\right]\}
\end{equation}
where the second bound holds only if $\Delta_r <
\frac{\sqrt{5}-1}{4}$.}

\subsection{Proof of Proposition 1}

\vskip.1in\noindent{\bf Proof :} For $n$ i.i.d. observations $X_i,
i=1,\ldots,n$, the KL discrepancy of the data is just $n$ times the
KL discrepancy for a single observation. Therefore, w.l.o.g. take $n
= 1$. Direct computation yields
\begin{equation}\label{eq:Sigma_inverse}
\Sigma^{-1} = (I - \sum_{\nu=1}^M \eta(\lambda_\nu) \theta_\nu
\theta_\nu^T).
\end{equation}
Hence, the log-likelihood function for a single observation is given
by
\begin{eqnarray}\label{eq:multi_log_likelihood}
\log f(x|\theta) &=& -\frac{N}{2}\log(2\pi) -
\frac{1}{2}\log|\Sigma|
- \frac{1}{2} x^T\Sigma^{-1}x \nonumber\\
&=& -\frac{N}{2}\log(2\pi) - \frac{1}{2}\sum_{\nu = 1}^M
\log(1+\lambda_\nu) - \frac{1}{2} (\langle x,x \rangle -
\sum_{\nu=1}^M \eta(\lambda_\nu) \langle x,\theta_\nu \rangle^2).
\end{eqnarray}
Recall that, if distributions $F_1$ and $F_2$ have density functions
$f_1$ and $f_2$, respectively, such that the support of $f_1$ is
contained in the support of $f_2$, then the Kullback-Leibler
discrepancy of $F_2$ from $F_1$, to be denoted by $K(F_1,F_2)$, is
given by
\begin{equation}\label{eq:KL_general}
K(F_1,F_2) = \int \log \frac{f_1(y)}{f_2(y)} f_1(y) dy.
\end{equation}
Hence, from (\ref{eq:multi_log_likelihood}),
\begin{eqnarray*}
K_{1,2} &=& \mathbb{E}_{\theta^{(1)}} (\log f(X|\theta^{(1)} - \log
f(X|\theta^{(2)}) \nonumber\\
&=& \frac{1}{2}  \sum_{\nu=1}^M \eta(\lambda_\nu)
[\mathbb{E}_{\theta^{(1)}}(\langle X , \theta_\nu^{(1)} \rangle)^2 -
\mathbb{E}_{\theta^{(1)}}(\langle X ,
\theta_\nu^{(2)} \rangle)^2] \nonumber\\
&=& \frac{1}{2}  \sum_{\nu=1}^M \eta(\lambda_\nu) [\langle
\theta_\nu^{(1)}, \Sigma_{(1)} \theta_\nu^{(1)}\rangle - \langle
\theta_\nu^{(2)}, \Sigma_{(1)} \theta_\nu^{(2)}\rangle] \nonumber\\
&=& \frac{1}{2}  \sum_{\nu=1}^M \eta(\lambda_\nu) [(\parallel
\theta_\nu^{(1)}\parallel^2 - \parallel
\theta_\nu^{(2)}\parallel^2)^2 + \sum_{\nu'=1}^M \lambda_{\nu'}
\{\parallel \theta_{\nu'}^{(1)}\parallel^2 - (\langle
\theta_{\nu'}^{(1)}, \theta_\nu^{(2)}\rangle)^2\}],
\end{eqnarray*}
which equals the RHS of (\ref{eq:multi_KL_div}), since the columns
of $\theta^{(j)}$ are orthonormal for each $j=1,2$.

\subsection{A counting lemma}\label{subsec:highdpca-counting}

{\bf Lemma 7:}\label{lemma:maximal_N_m_set} \textit{Suppose that
$m$, $N$ are positive integers, such that $m \to \infty$ as $N \to
\infty$ and $m = o(N)$. Let $\widetilde Z$ be the maximal set of
points in $\mathbb{R}^N$ satisfying the following conditions:
\begin{itemize}
\item[(i)] for each $\mathbf{z} = (z_1,\ldots,z_N) \in \widetilde Z$,
$z_i \in \{0,1\}$
  for all $i=1,\ldots,N$,
\item[(ii)] for each $\mathbf{z}  \in \widetilde Z$,
exactly $m$ of coordinates of
  $\mathbf{z}$ are 1,
\item[(iii)] for every pair
${\bf z}$ and ${\bf z}'$ in $\widetilde Z$, $z_i = z_i'$ for at
  most $\left[\frac{m_0}{2}\right] =: k(m_0) - 1$ (i.e. $k(m_0)$ is the largest
  integer $\leq m_0/2+1$)
nonzero coordinates, where $m_0 = [\beta m]$, for some $\beta \in
(0,1)$.
\end{itemize}
Then cardinality of $\widetilde Z$ is at least $\exp([N{\cal
E}(\frac{\beta m}{2N})- 2m{\cal E}(\frac{\beta}{2})](1+o(1)))$ where
${\cal E}(x)$ is the Shannon entropy function.}

\vskip.1in\noindent{\bf Proof : } Trivially, $\widetilde Z \subset
Z^*$, where $Z^*$ is the set of {\it all} points $\mathbf{z}$
satisfying (i) and (ii). Thus, $|\widetilde Z| < | Z^*| = {N \choose
m}$. On the other hand, for every point $\mathbf{z} \in Z^*$ there
are \textit{at most}
$$
g(N,m,m_0) = {m \choose k(m_0)} {{N - k(m_0)}\choose {m - k(m_0)}}
$$
points $\mathbf{w} \in Z^*$ such that \textit{at least} $k(m_0)$
nonzero coordinates of $\mathbf{z}$ and $\mathbf{w}$ match. This is
because, one can fix the $m$ nonzero coordinates of $\mathbf{z}$ and
demand that in $k(m_0)$ of those coordinates $w_i$ must equal 1.
Other $m-k(m_0)$ nonzero coordinates of $\mathbf{w}$ can therefore
be chosen from the rest $N - k(m_0)$ coordinates. Then, by the
maximality of $\widetilde Z$, as $N \to \infty$,
\begin{eqnarray}\label{eq:card_Z_*}
|\widetilde Z | &\geq & {N \choose m} g(N,m,m_0)^{-1} \nonumber\\
&=& \frac{N!}{(N- m)! m!} ~\frac{k(m_0)!(m - k(m_0))!}{m!}
~\frac{(m - k(m_0))!(N-m)!}{(N-k(m_0))!} \nonumber\\
&=& \frac{N!}{k(m_0)!(N - k(m_0))!}
~\left(\frac{k(m_0)!(m-k(m_0))!}{m!}\right)^2 \nonumber\\
&\sim& \sqrt{2\pi}\sqrt{k(m_0)} ~
\frac{(m-k(m_0))N^{1/2}}{m(N-k(m_0))^{1/2}}~
(\frac{N}{k(m_0)})^{k(m_0)}
(\frac{N}{N - k(m_0)})^{N-k(m_0)} \nonumber\\
&& ~~ \times ~
[(\frac{k(m_0)}{m})^{k(m_0)}(\frac{m-k(m_0)}{m})^{m-k(m_0)}]^2
~~~(\mbox{by
  Stirling's formula})
\nonumber\\
&=& \sqrt{2\pi}\sqrt{\frac{\beta m}{2}} \exp\left[N {\cal
E}\left({\frac{\beta m}{2N}}\right)(1+o(1))\right] \exp\left[- 2m
{\cal E}\left(\frac{\beta}{2}\right)(1+o(1))\right].
\end{eqnarray}
Where the last equality is because, for large $m$, $\frac{m_0}{m}
\sim \frac{\beta}{2}$.

\subsection{Some auxiliary lemmas}\label{subsec:aspca-lemmas}

In the following lemmas we provide probabilistic bounds for the
deviations of certain quadratic forms that arise in the analysis of
the residual terms in the expansion of $\widehat \theta_\nu$. Many
of these involve the random sets, either $\widehat I_{1,n}$ or
$\widehat I_{2,n}$, of coordinates that are selected under the ASPCA
scheme. It will be assumed that the quantities involved are all
measurable w.r.t. the joint distribution of $\mathbf{Z}$ and
$v_1,\ldots,v_M$, though it will not be made explicit in the
description or the proof of the lemmas. The bounds hold uniformly in
$\theta \in \Theta_q^M(C_1,\ldots,C_M)$.

\vskip.15in\noindent{\bf Lemma 8:}\label{lemma:Z_k_Y_A_bound}
\textit{Let $\epsilon_n
> 0$. Let $A$ denote the random set $\widehat I_{n,1}$, and $A_- =
I_{1,n}^-$ and $A_+ = I_{1,n}^+$. Assume that $A_- \subset
A_+\setminus \{k\}$, for some $1\leq k \leq N$. For any subset $C$
of $\{1,\ldots,N\}$, let $Y_C := Y_C(\mathbf{Z}_C,V)$ be a random
vector jointly measurable w.r.t. $\mathbf{Z}_C$ and $V =
[v_1:\ldots:v_M]$. Assume that for each $C$, either
$\mathbb{P}_V(Y_C = 0) = 0$ a.e. $V$, or $\mathbb{P}_V(Y_C = 0) = 1$
a.e. $V$, where $\mathbb{P}_V$ denotes the conditional probability
w.r.t. $V$. Let $W_{k,C} = \langle Z_k, \frac{Y_C}{\parallel Y_C
\parallel} \rangle$ if $Y_C \neq 0$, and $W_{k,C} = 0$ otherwise.
Then,
\begin{equation}\label{eq:dep_inner_prod_bound}
\mathbb{P}\left( |W_{k,A}| > \epsilon_n , A_- \subset A \subset A_+
\setminus \{k\}, ~
\parallel V \parallel \leq \beta_n\right)
\leq \frac{2}{a_n} \Phi(-\epsilon_n),
\end{equation}
where $\beta_n$ is such that, on $\{\parallel V \parallel \leq
\beta_n\}$, a.e. $V$,
\begin{equation}\label{eq:dep_inner_prod_cond}
\mathbb{P}_V(\widehat \sigma_{kk} \leq 1+ \gamma_{1,n}) \geq a_n >
0.
\end{equation}}

\vskip.15in\noindent{\bf Lemma 9:}
\label{lemma:Z_A_V_mu_theta_nu_bound} \textit{Let $A$ be a random
subset of $\{1,\ldots,N\}$ and $A_- \subset A_+$ be two non-random
subsets of $\{1,\ldots,N\}$. Let, $k_{\pm}$ denote the size of the
set $A_{\pm}$, and
$$
\epsilon_n =  \sqrt{c_1 \log n} + \parallel \theta_{\nu,A_-^c}
\parallel \sqrt{c_1 \log n + 2 k_+ \log 2},
$$
for some $c_1 > 0$. Then, for all $1\leq \mu \leq M$,
\begin{equation}\label{eq:Z_A_v_theta_A_bound}
\mathbb{P}\left(|\langle \frac{\mathbf{Z}_A v_\mu}{\parallel v_\mu
    \parallel}, \theta_{\nu,A}\rangle| > \epsilon_n, A_- \subset A
  \subset A_+\right)
\leq \frac{4  n^{-c_1/2}}{\sqrt{2\pi}{\sqrt{c_1 \log n}}} .
\end{equation}}

\vskip.15in\noindent{\bf Lemma 10:}
\label{lemma:Z_A_theta_nu_Z_A_theta_mu_bound} \textit{Let $A$,
$A_\pm$, $k_\pm$, and $\overline{A}_-$ be as in Lemma 9.
Let
$$
\epsilon_n = \parallel \theta_{\mu,A_-^c}\parallel
(1+\sqrt{\frac{k_+}{n}}+\sqrt{\frac{c_2 \log n}{n}}) \sqrt{\frac{c_1
\log n + 2k_+\log 2}{n}},
$$
where $c_1, c_2 > 0$. Then
\begin{equation}\label{eq:Z_A_theta_mu_nu_bound}
\mathbb{P}\left( |\frac{1}{n} \langle \mathbf{Z}_{A_-}^T
  \theta_{\nu,A_-}, \mathbf{Z}_{\overline{A}_-}^T
  \theta_{\mu,\overline{A}_-} \rangle | > \epsilon_n, ~A_- \subset A
  \subset A\right)
\leq \frac{2n^{-c_1/2}}{\sqrt{2\pi}\sqrt{c_1 \log
    n}} + n^{-c_2/2}.
\end{equation}}

\vskip.15in\noindent{\bf Lemma 11:}
\label{lemma:Z_A_theta_nu_A_quad_bound} \textit{Let $A$, $A_\pm$,
$k_\pm$ be as in Lemma 9.
Let,
$t_n = 6(\frac{k_+}{n}\vee 1)\sqrt{\frac {\log(n \vee k_+)}{n\vee
k_+}}$. Let
\begin{eqnarray*}
\epsilon_n &=& \sqrt{\frac{c_1 \log n}{n}} +
\parallel \theta_{\nu,A_-^c} \parallel^2 (2\sqrt{\frac{k_+}{n}}
+ \frac{k_+}{n} + \sqrt{c_2/2} t_n) \\
&& ~~~+ 2 \parallel \theta_{\nu,A_-^c}\parallel (1 +
\sqrt{\frac{k_+}{n}} + \sqrt{\frac{c_2 \log
    n}{n}}) \sqrt{\frac{c_2 \log n + 2k_+\log 2}{n}},
\end{eqnarray*}
for some $c_1, c_2 > 0$. Then there is an $n(c_2) \geq 16$ such
that, for $n \geq n(c_2)$, $\sqrt{c_2/2} t_n < 1/2$, and
\begin{eqnarray}\label{eq:Z_A_theta_A_norm_bound}
&&\mathbb{P}\left(|\frac{1}{n} \parallel \mathbf{Z}_A^T
\theta_{\nu,A}
\parallel^2 - \parallel \theta_{\nu,A} \parallel^2 | > \epsilon_n,
~ A_- \subset A
\subset A_+ \right) \nonumber\\
&\leq& 2n^{-c_1/4} + \frac{2n^{-c_2/2}}{\sqrt{2\pi}\sqrt{c_1 \log
    n}} + n^{-c_2/2} + 2(n \vee k_+)^{-c_2/2}.
\end{eqnarray}}

\vskip.15in\noindent{\bf Lemma 12:}
\label{lemma:Z_A_theta_nu_Z_A_theta_mu_centered} \textit{Let $A$,
$A_\pm$, $k_\pm$ be as in Lemma 9.
Let, $\mu \neq \nu$, and for
some $t > 0$,
\begin{eqnarray*}
\epsilon_n &=& \sqrt{\frac{c_1 \log n}{n}} + (\parallel
\theta_{\mu,A_-^c}\parallel + \parallel \theta_{\nu,A_-^c}\parallel)
(1 + \sqrt{\frac{k_+}{n}}+\sqrt{\frac{c_3\log n}{n}})
\sqrt{\frac{c_2 \log n + 2k_+\log 2}{n}} \\
&&+ \parallel \theta_{\nu,A_-^c} \parallel \parallel
\theta_{\mu,A_-^c} \parallel \sqrt{\frac{c_2 \log n}{n}} +
\parallel \theta_{\nu,A_-^c} \parallel \parallel \theta_{\mu,A_-^c} \parallel
(2\sqrt{\frac{k_+}{n}} + \frac{k_+}{n}+\sqrt{c_3/2}t_n),
\end{eqnarray*}
where $c_1,c_2,c_3 > 0$ and $t_n$ is as in Lemma 11.
Then, there is $n(c_3) \geq
16$ such that, for $n \geq n(c_3)$, $\sqrt{c_3}t_n < \frac{1}{2}$,
and
\begin{eqnarray}\label{eq:Z_A_theta_mu_nu_A_bound}
&&\mathbb{P}\left(|\frac{1}{n} \langle \mathbf{Z}_A^T
\theta_{\nu,A}, \mathbf{Z}_A^T \theta_{\mu,A} \rangle - \langle
\theta_{\nu,A}, \theta_{\mu,A} \rangle | > \epsilon_n, A_- \subset A
\subset A_+ \right) \nonumber\\
&\leq& 2n^{-3c_1/2+O(\frac{\log n}{n})} + 2n^{-c_2/4} +
\frac{2n^{-c_2/2}}{\sqrt{2\pi}\sqrt{c_1 \log
     n}} + n^{-c_3/2} + 2(n \vee k_+)^{-c_3/2}.
\end{eqnarray}}

\vskip.1in\noindent{\bf Lemma 13:}
\label{lemma:Z_A_Z_A_bar_theta_mu_bound} \textit{Let $A$, $A_\pm$,
$k_\pm$ be as in Lemma 9.
Let,
\begin{eqnarray*}
\epsilon_n &=& 2\parallel \theta_{\mu,A_-^c}\parallel(1 +
\sqrt{\frac{k_+}{n}}+\sqrt{\frac{c_2\log n}{n}})
\sqrt{\frac{k_+}{n}} \left(1+ \sqrt{\log 2 + \frac{c_1\log
n}{4k_+}}\right)^{1/2} ,
\end{eqnarray*}
where $c_1,c_2 > 0$. Also, suppose that $k_+ \geq 16$. Then,
\begin{equation}\label{eq:Z_A_Z_A_bar_theta_mu_bound}
\mathbb{P}(\frac{1}{n}\parallel \mathbf{Z}_{A_-}
\mathbf{Z}_{\overline{A}_-}^T \theta_{\mu,\overline{A}_-}\parallel >
\epsilon_n, A_- \subset A \subset A_+) \leq n^{-c_1/4} + n^{-c_2/2}.
\end{equation}}

\subsection{Deviation of quadratic
forms}\label{subsec:deviation_quad_form}

The following lemma is due to Johnstone (2001b).

\vskip.15in\noindent{\bf Lemma 14:}
\label{lemma:chi_square_large_dev} \textit{Let $\chi_{(n)}^2$ denote
a Chi-square random variable with $n$ degrees of freedom. Then,
\begin{eqnarray}
 \mathbb{P}(\chi_{(n)}^2 > n(1+\epsilon) ) &\leq& e^{-3n\epsilon^2/16}
\qquad (0 < \epsilon < \frac{1}{2}),\label{eq:large_dev_chisq_1}\\
\mathbb{P}(\chi_{(n)}^2 < n(1-\epsilon) ) &\leq& e^{-n\epsilon^2/4}
\qquad (0 < \epsilon < 1),\label{eq:large_dev_chisq_2}\\
\mathbb{P}(\chi_{(n)}^2 > n(1+\epsilon) ) &\leq&
\frac{\sqrt{2}}{\epsilon
  \sqrt{n}}e^{-n\epsilon^2/4}
\qquad (0 < \epsilon < n^{1/16}, n \geq 16).
\label{eq:large_dev_chisq_3}
\end{eqnarray}}

\noindent The following lemma is from Johnstone and Lu (2004).
\vskip.15in\noindent{\bf Lemma 15:}
\label{lemma:quad_form_large_dev} \textit{Let
  $y_{1i},y_{2i},i=1,\ldots,n$ be two sequences of mutually
  independent, i.i.d. $N(0,1)$ random variables. Then for large $n$
  and any $b$ s.t. $0 < b \ll \sqrt{n}$,
\begin{equation}\label{eq:large_dev_covar}
\mathbb{P}(|\frac{1}{n} \sum_{i=1}^n y_{1i} y_{2i}|>\sqrt{b/n}) \leq
2\exp\{-\frac{3b}{2} + O(n^{-1}b^2)\}.
\end{equation}}

\section*{Reference}

\begin{enumerate}
\item
Anderson, T. W. (1963) :
     Asymptotic theory of principal component analysis,
     \textit{Annals of Mathematical Statistics}, {\bf 34},
     122-148.
\item
Bai, J. (2003) :
     Factor models for large dimensions,
     \textit{Econometrica}, {\bf 71}, 135-171.
\item
Bair, E., Hastie, T., Paul, D. and Tibshirani, R. (2006) :
     Prediction by supervised principal components,
     \textit{Journal of the American Statistical Association},
     {\bf 101}, 119-137.
\item
Birg\'{e}, L. (2001) :
     A new look at an old result : Fano's lemma,
     \textit{Technical Report}, Universit\'{e} Paris 6.
\item
Boente, G. and Fraiman, R. (2000) :
     Kernel-based functional principal components,
     \textit{Statistics and Probability Letters}, {\bf 48}, 335-345.
\item
Buja, A. and  Hastie, T. and Tibshirani, R. (1995) :
     Penalized discriminant analysis,
     \textit{Annals of Statistics}, {\bf 23}, 73-102.
\item
Cassou, C., Deser, C., Terraty, L., Hurrell, J. W. and
Dr\'{e}villon, M. (2004) :
     Summer sea surface temperature conditions in the
North Atlantic and their impact upon the atmospheric circulation in
early winter,
     \textit{Journal of Climate}, {\bf 17}, 3349-3363.
\item
Cardot, H. (2000) :
     Nonparametric estimation of smoothed principal components
analysis of sampled noisy functions,
     \textit{Journal of Nonparametric Statistics}, {\bf 12},
     503-538.
\item
Cardot, H., Ferraty, F. and Sarda, P. (2003) :
     Spline estimators for the functional linear model,
     \textit{Statistica Sinica}, {\bf 13}, 571-591.
\item
Chiou, J.-M., M\"{u}ller, H.-G. and Wang, J.-L. (2004) :
     Functional response model, \textit{Statistica Sinica}, {\bf 14},
     675-693.
\item
Cootes, T. F., Edwards, G. J. and Taylor, C. J. (2001) :
     Active appearance models, \textit{IEEE Transactions on Pattern Analysis and Machine
Intelligence}, {\bf 23}, 681-685.
\item
Corti, S., Molteni, F. and Palmer, T. N. (1999) :
     Signature of recent climate change in frequencies of natural atmospheric circulation regimes,
     \textit{Nature}, {\bf 398}, 799-802.
\item
Davidson, K. R. and Szarek, S. (2001) :
     Local operator theory, random matrices and Banach spaces,
     in \textit{Handbook on the Geometry of Banach Spaces}, {\bf 1},
     \textit{Eds. Johnson, W. B. and Lendenstrauss, J.}, 317-366, Elsevier
     Science.
\item
Dey, D. K. and Srinivasan, C. (1985) :
     Estimation of a covariance matrix under Stein's loss,
     \textit{Annals of Statistics}, {\bf 13}, 1581-1591.
\item
Donoho, D. L. (1993) :
     Unconditional bases are optimal bases for data compression
and statistical estimation,
     \textit{Applied and Computational Harmonic Analysis}, {\bf 1}, 100-115.
\item
Eaton, M. L. and Tyler, D. E. (1991) :
     On Wielandt's inequality and its application to the
asymptotic distribution of a random symmetric matrix,
     \textit{Annals of Statistics}, {\bf 19}, 260-271.
\item
Efron, B. and Morris, C. (1976) :
     Multivariate empirical Bayes estimation of covariance matrices,
     \textit{Annals of Statistics}, {\bf 4}, 22-32.
\item
Haff, L. R. (1980) :
     Empirical Bayes estimation of the multivariate normal
covariance matrix, \textit{Annals of Statistics}, {\bf 8},
     586-597.
\item
Hall, P. (1992) :
     \textit{The Bootstrap and Edgeworth Expansion},
     Springer-Verlag.
\item
Hall, P. and Horowitz, J. L. (2004) :
     Methodology and convergence rates for functional linear regression,
     \textit{Manuscript}.
\item
Hall, P. and Hosseini-Nasab, M. (2006) :
     On properties of functional
     principal components analysis, \textit{Journal of Royal Statistical
     Society, Series B}, {\bf 68}, 109-125.
\item
Tony Cai, T. and Hall, P. (2005) :
     Prediction in functional linear regression,
     \textit{Manuscript}.
\item
Hoyle, D. and Rattray, M. (2003) :
     Limiting form of the sample covariance eigenspectrum in PCA and kernel
     PCA,
     \textit{Advances in Neural Information Processing Systems},
     {\bf 16}.
\item
Hoyle, D. and Rattray, M. (2004) :
     Principal component analysis eigenvalue spectra from data with symmetry
     breaking structure, \textit{Physical Review E}, {\bf 69},
     026124.
\item
Johnstone, I. M. (2001) :
     On the distribution of the largest principal component,
     \textit{Annals of Statistics}, {\bf 29}, 295-327.
\item
Johnstone, I. M. (2001b) :
     Chi square oracle inequalities, in
     \textit{Festschrift for William R. van Zwet}, {\bf 36},
     \textit{Eds. de Gunst, M., Klaassen, C. and Waart, A. van der},
     399-418, Institute of Mathematical Statistics.
\item
Johnstone, I. M. (2002) :
     \textit{Function estimation and gaussian sequence models},
     Book Manuscript.
\item
Johnstone, I. M. and Lu, A. Y. (2004) :
     Sparse principal component analysis,
     \textit{Technical Report}, Stanford University.
\item
Kato, T. (1980) :
     \textit{Perturbation Theory of Linear Operators},
     Springer-Verlag.
\item
Kneip, A. (1994) :
     Nonparametric estimation of common regressors for similar curve data,
     \textit{Annals of Statistics}, {\bf 22}, 1386-1427.
\item
Kneip, A. and Utikal, K. J. (2001) :
     Inference for density families using functional principal component analysis,
     \textit{Journal of the American Statistical Association}, {\bf 96},
     519-542.
\item
Laloux, L., Cizeau, P., Bouchaud, J. P. and Potters, M. (2000) :
     Random matrix theory and financial correlations,
     \textit{International Journal of Theoretical and Applied Finance},
     {\bf 3}.
\item
Loh, W.-L. (1988) :
     Estimating covariance matrices,
     \textit{Ph. D. Thesis}, Stanford University.
\item
Lu, A. Y. (2002) :
     Sparse principal component analysis for functional data,
     \textit{Ph. D. Thesis}, Stanford University.
\item
Muirhead, R. J. (1982) :
     \textit{Aspects of Multivariate Statistical Theory},
     John Wiley \& Sons, Inc.
\item
Paul, D. (2005) : Nonparametric estimation of principal components,
     \textit{Ph. D. Thesis}, Stanford University.
\item
Paul, D. and Johnstone, I. M. (2004) :
     Estimation of principal components through coordinate selection,
     \textit{Technical Report}, Stanford University.
\item
Preisendorfer, R. W. (1988) :
     \textit{Principal component analysis in meteorology and oceanography},
     Elsevier, New York.
\item
Ramsay, J. O.  and Silverman, B. W. (1997) :
     \textit{Functional Data Analysis}, Springer-Verlag.
\item
Ramsay, J. O.  and Silverman, B. W. (2002) :
     \textit{Applied Functional Data Analysis : Methods and Case Studies},
     Springer-Verlag.
\item
Spellman, P.T., Sherlock, G., Zhang, M. Q.,  Iyer, V. R.,  Anders,
K., Eisen, M. B., Brown, P. O., Botstein, D. and Futcher, B. (1998)
:
     Comprehensive identification of cell cycle-regulated genes of the
     yeast saccharomyces cerevisiae by microarray hybridization,
     \textit{Molecular Biology of the Cell}, {\bf 9}, 3273-3297.
\item
Stegmann, M. B. and Gomez, D. D. (2002) :
     A brief introduction to statistical shape analysis,
     \textit{Lecture notes}, Technical University of Denmark.
\item
Telatar, E. (1999) :
     Capacity of multi-antenna Gaussian channels,
     \textit{European Transactions on Telecommunications}, {\bf 10},
     585-595.
\item
Tulino, A. M. and Verdu, S. (2004) :
     Random matrices and wireless communications,
     \textit{Foundations and Trends in Communications and Information Theory},
     {\bf 1}.
\item
Tyler, D. E. (1983) :
     The asymptotic distribution of principal component roots under
local alternatives to multiple roots,
     \textit{Annals of Statistics}, {\bf 11}, 1232-1242.
\item
Vogt, F., Dable, B., Cramer, J. and Booksh, K. (2004) :
     Recent advancements in chemometrics for smart sensors,
     \textit{The Analyst}, {\bf 129}, 492-502.
\item
Wickerhauser, M. V. (1994) :
     \textit{Adapted Wavelet Analysis from Theory to Software},
     A K Peters, Ltd.
\item
Yang, Y. and Barron, A. (1999) :
     Information-theoretic determination of minimax rates of convergence,
     \textit{Annals of Statistics}, {\bf 27}, 1564-1599.
\item
Zhao, X., Marron, J. S. and Wells, M. T. (2004) :
     The functional data analysis view of longitudinal data,
     \textit{Statistica Sinica}, {\bf 14}, 789-808.
\item
Zong, C. (1999) :
     \textit{Sphere Packings}, Springer-Verlag.

\end{enumerate}

\end{document}